\documentclass{amsart}

\usepackage{amsmath}

\usepackage{epsfig}
\usepackage{graphicx}
\usepackage{amsbsy}
\usepackage{latexsym}
\usepackage[mathscr]{eucal}
\usepackage{amsfonts,amsmath,amssymb}
\usepackage{hyperref}
\usepackage{algorithm}
\usepackage{algpseudocode}

\usepackage[latin1]{inputenc}
\usepackage{color}
\usepackage{pgfplots}
\usepackage{subfig}
\newtheorem{lemma}{Lemma}[section]
\newtheorem{prop}[lemma]{Proposition}

\newtheorem{theorem}[lemma]{Theorem}
\newtheorem{rem}[lemma]{Remark}

\newcommand{\bo}{B}
\newcommand{\e}{\varepsilon}
\newcommand{\mup}{y}
\newcommand{\N}{\mathbb{N}}
\newcommand{\R}{\mathbb{R}}
\newcommand{\cK}{{\mathcal{K}}}
\newcommand{\cM}{{\mathcal{M}}}
\newcommand{\cN}{{\mathcal{N}}}
\newcommand{\cY}{\mathcal{Y}}

\newcommand{\dist}{\mathop{\rm max\,dist}}
\newcommand{\Span}{\mathop{\rm span}}

\DeclareMathOperator*{\argmin}{argmin}
\DeclareMathOperator*{\argmax}{argmax}

\newcommand{\SGA}{\textsc{SGA}}
\newcommand{\WGA}{\textsc{WGA}}
\newcommand{\DoubGA}{\textsc{SGA-dou}}

\makeatletter
\@addtoreset{equation}{section}
\makeatother

\begin{document}
\title{How To Best Sample a Solution Manifold?}
\thanks{This work has been supported in part by the DFG SFB-Transregio 40,
and by the DFG Research Group 1779, the Excellence Initiative of the German Federal and State Governments,
and NSF grant DMS 1222390}
\author{Wolfgang Dahmen}\address{Institut f\"ur Geometrie und Praktische Mathematik, RWTH Aachen, Germany, e-mail: dahmen@igpm.rwth-aachen.de}
\date{}
\begin{abstract}
Model reduction attempts to guarantee a desired ``model quality'', e.g. given in terms of
accuracy requirements, with as small a model size as possible. This article highlights some recent
developments concerning this issue for the so called Reduced Basis Method (RBM) for models based
on parameter dependent families of PDEs. In this context the key task is to sample the {\em solution manifold}
at judiciously chosen parameter values usually determined in a {\em greedy fashion}.
The corresponding {\em space growth} concepts are closely related to so called {\em weak greedy} algorithms in Hilbert and Banach spaces
which can be shown to give rise to convergence rates comparable to the best possible rates, namely the {\em Kolmogorov $n$-widths} rates.
Such algorithms can be interpreted as {\em adaptive sampling} strategies for approximating compact sets in Hilbert spaces.
We briefly discuss the results most relevant for the present RBM context.
The applicability of the results for weak greedy algorithms has however been confined so far essentially to
well-conditioned coercive problems. A critical issue is therefore an extension of these concepts to a wider range of problem classes for which the conventional methods
do not work well. A second main topic of this article is therefore to outline recent developments of RBMs that do realize $n$-width rates
for a much wider class of variational problems covering indefinite or singularly perturbed unsymmetric problems.
A key element in this context is the design of {\em well-conditioned variational formulations}
and their numerical treatment via saddle point formulations. We conclude with some remarks concerning the relevance
of uniformly approximating the whole solution manifold also when the {\em quantity of interest} is only of a {\em functional} of the parameter dependent solutions.
\end{abstract}
\subjclass{65J10, 65N12, 65N15, 35B30}
\keywords{Tight surrogates, stable variational formulations, saddle point problems, double greedy schemes, greedy stabilization,
rate-optimality, transport equations, convection-diffusion equations.}
\maketitle

%%%%%%%%%%%%%%%%%%%%%%%%%%%%%
\section{Introduction}\label{intro}
%%%%%%%%%%%%%%%%%%%%%%%%%%%%%
Many engineering applications revolve around the task of identifying a configuration that in some sense best fits certain objective criteria
under certain constraints. Such design or optimization problems typically involve (sometimes many)
{\em parameters} that need to be chosen so as to satisfy given
optimality criteria. An optimization over such a parameter domain usually requires a frequent evaluation of
the states under consideration which typically means to frequently solve a {\em parameter dependent} family of operator equations
\begin{equation}
\label{1.1}
B_\mup u = f, \quad \mup\in\cY.
\end{equation}
In what follows the parameter set $\cY$ is always assumed to be a compact subset of $\R^p$ for some fixed $p\in \N$ and
$B_\mup$ should be thought of as a (linear) partial differential operator whose coefficients depend on the parameters $\mup\in \cY$.
Moreover, $B_\mup$ is viewed as an operator taking some Hilbert space $U$ one-to-one and onto the {\em normed dual} $V'$ of some (appropriate)
Hilbert space $V$ where $U$ and $V$ are identified through a variational formulation of \eqref{1.1} as detailed later, see for instance \eqref{condiff}.
Recall also that the normed dual $V'$ is endowed with the norm
\begin{equation}
\label{dualnorm}
\|w\|_{V'}:= \sup_{v\in V, v\neq 0}\frac{\langle w,v\rangle}{\|v\|_V},
\end{equation}
where $\langle \cdot,\cdot\rangle$ denotes the dual pairing between $V$ and $V'$.

Given a parametric model \eqref{1.1} the above mentioned design or optimization problems concern now
the {\em states} $u(\mup)\in U$ which, as a function of the parameters $\mup\in \cY$, form what we refer to as the {\em solution manifold}
\begin{equation}
\label{sol-man}
\cM := \{u(\mup):= B_\mup^{-1}f :\mup\in\cY\}.
\end{equation}
Examples of \eqref{1.1} arise, for instance, in geometry optimization when a transformation of a variable finitely parametrized domain to a reference domain introduces parameter dependent coefficients of the underyling partial differential equation (PDE) over such domains, see e.g. \cite{GV}. Parameters could describe conductivity, viscosity or convection directions, see e.g. \cite{DPW,RHP,SVHDNP}. As an extreme case,
parametrizing the random diffusion coefficients in a stochastic PDE e.g., by Karhunen-Loew or polynomial chaos expansions,
leads to a deterministic parametric PDE involving, in principle, even {\em infinitely} many parameters, $p=\infty$, see e.g. \cite{CDS1} and the
literature cited there.
We will, however, not treat this particular problem class here any further since, as will be explained later, it poses different conceptual obstructions than
those in the focus of this paper, namely the absence of ellipticity which makes
conventional strategies fail.
In particular, we shall explain why for other relevant problem classes,
e.g. those dominated by transport processes, $\cM$ is not ``as visible'' as for elliptic problems and how to restore ``full visibility''.

\subsection{General Context - Reduced Basis Method}
A conventional way of searching for a specific state in $\cM$ or optimize over $\cM$ is to compute approximate solutions of \eqref{1.1}
possibly for a large number of parameters $\mup$. Such approximations would then reside in a sufficiently large trial space $U_\cN\subset U$
of dimension $\cN$, typically a finite element space. Ideally one would try to assure that $U_\cN$ is large enough to warrant sufficient accuracy
of whatever conclusions are to be drawn from such a discretization. A common terminology in reduced order modeling refers to $U_\cN$
as the {\em truth space} providing accurate computable information. Of course, each such parameter query in $U_\cN$ is a computationally
expensive task so that many such queries, especially in an online context, would be practically infeasible. On the other hand,
solving for each $y\in \cY$ a problem in $U_\cN$
would just treat each solution $u(\mup)$ as some ``point'' in the infinite dimensional space $U$,
viz. in the very large finite dimensional space $U_\cN$. This
disregards the fact that all these points actually belong to a possibly much thinner and more coherent set, namely the low dimensional manifold $\cM$ which, for compact $\cY$ and well posed problems
\eqref{1.1}, is compact. Moreover, if the solutions $u(\mup)$, as functions of $\mup\in \cY$, depend smoothly on $\mup$ there is hope that one
can approximate all elements of $\cM$ uniformly over $\cY$ with respect to the Hilbert space norm $\|\cdot\|_U$ by a
relatively small but judiceously chosen linear space $U_n$. Here ``small'' means that $n={\rm dim}\,U_n$ is significantly smaller than $\cN={\rm dim}\,U_\cN$, often
by orders of magnitude. As detailed later the classical notion of {\em Kolmogorov $n$-widths} quantifies how well a compact set in a Banach space can be approximated in the corresponding Banach norm by a linear space and therefore can be used as a {\em benchmark} for
the effectiveness of a model reduction strategy.

Specifically, the core objective of the {\em Reduced Basis Method} (RBM) is to find for a given {\em target accuracy} $\e$ a possibly small
number $n=n(\e)$ of basis functions
$\phi_j, j=0,\ldots,n,$ whose linear combinations approximate each $u\in \cM$ within accuracy at least $\e$.
This means
that ideally for each $\mup\in \cY$ one can find coefficients $c_j(\mup)$ such that the expansion
\begin{equation}
\label{sep}
u_n(x,\mup):= \sum_{j=0}^{n(\e)}c_j(\mup)\phi_j(x)
\end{equation}
satisfies
\begin{equation}
\label{target}
\|u(\mup) - u_n(\mup)\|_U \le \e,\quad \mup\in \cY.
\end{equation}
Thus, projecting \eqref{1.1} into the small
space $U_n := {\rm span}\, \{\phi_0,\ldots,\phi_n\}$ reduces each parameter query to solving a small $n\times n$ system of equations
where typically $n\ll \cN$.

\subsection{Goal Orientation}
Recall that the actual goal of reduced modeling is often not to recover the full fields $u(\mup)\in \cM$ but only some {\em quantity of interest} $I(\mup)$
typically given as a {\em functional} $I(\mup):= \ell(u(\mup))$ of $u(\mup)$ where $\ell \in U'$.
Asking just the value of such a functional is possibly a weaker request than approximating all of $u(\mup)$ in the norm $\|\cdot\|_U$.
In other words, one may have $|\ell(u(\mup))- \ell(u_n(\mup))|\le \e$ without insisting on the validity of \eqref{target} for a tolerance of
roughly the same size.
Of course, one would like to
exploit this in favor of online efficiency. Duality methods as used in the context of {\em goal-oriented} finite element methods \cite{BR}
are indeed known to offer ways of economizing the approximate evaluation of functionals. Such concepts apply in the RBM context as well,
see e.g. \cite{Prud,Grepl1}. However, as we shall point out later, guaranteeing that $|\ell(u(\mup))- \ell(u_n(\mup))|\le \e$ holds for $\mup\in \cY$,
ultimately reduces to tasks of the type \eqref{target} as well. So, in summary, understanding how to ensure \eqref{target} for possibly small $n(\e)$
remains the core issue and therefore guides the subsequent discussions.

Postponing for a moment the issue of how to actually compute the $\phi_j$, it is clear that they should intrinsically depend on $\cM$ rendering the whole process
highly nonlinear.
To put the above approach first into perspective, viewing $u(x,\mup)$ as
a function of the spatial variables $x$ and of the parameters $\mup$, \eqref{sep} is just {\em separation of variables} where the factors
$c_j(\mup)$, $\phi_j(x)$ are a priori unknown. It is perhaps worth stressing though that, in contrast to other attempts to find
good {\em tensor approximation}, in the RBM context explicit representations are only computed for the spatial factors $\phi_j$
while for each $\mup$ the weight $c_j(\mup)$ has to be {\em computed} by solving a small system in the reduced space $U_n$.
Thus, the computation of $\{\phi_0,\ldots,\phi_{n(\e)}\}$ could be interpreted as {\em dictionary learning} and, loosely speaking, $n=n(\e)$ being relatively small for a given
target accuracy, means that all elements in $\cM$ are {\em approximately sparse} with respect to the dictionary $\{\phi_0,\ldots,\phi_n, \ldots\}$.

The methodology just outlined has been pioneered by Y. Maday, T.A. Patera and collaborators,
see e.g. \cite{Buffaetal,{Prud},{SVHDNP},RHP}. As indicated before, RBM is one variant of a {\em model oder reduction} paradigm that is specially tailored to parameter dependent
problems. Among its distinguishing constituents one can name the following. There is usually a careful division of the overall
computational work into an {\em offline phase}, which could be computationally intense but should remain managable, and
an {\em online phase} which should be executable with highest efficiency taking advantage of a
precomputed basis and matrix assemblations during the offline phase. It is important to note that
while the offline phase is accepted to be computationally expensive it should remain {\em offline-feasible} in the sense
that a possibly extensive search
over the parameter domain $\cY$ in the offline phase requires for each query solving only problems in the small reduced space.
Under which circumstances this is possible and how to realize such division concepts has been worked out
in the literature, see e.g. \cite{SVHDNP,RHP}. Here we are content with stressing that an important role is played by
the way how the operator $B_\mup$ depends on the parameter $\mup$, namely in an {\em affine} way as stated in \eqref{affine}
later below.
Second, and this is perhaps the most distinguishing constituent, along with each solution in the reduced model one strives to provide
a {\em certificate} of accuracy, i.e., computed bounds for incurred error tolerances \cite{SVHDNP,RHP}.

\subsection{Central Objectives}
When trying to quantify the performance of such methods aside from the above mentioned structural and data organization aspects, among others, the following questions come to mind:

(i) for which type
of problems do such methods work very well in the sense that $n(\e)$ in \eqref{target} grows only slowly when $\e$ decreases? This
concerns quantifying the sparsity of solutions.

(ii) How can one compute reduced bases $\{\phi_0,\ldots,\phi_{n(\e)}\}$ for which $n(\e)$ is {\em nearly minimal} in a sense to
be made precise below?

Of course, the better the sparsity quantified by (i) the better could be the pay-off of an RBM.
However, as one my expect, an answer to (i) depends strongly on
the problem under consideration. This is illustrated also by the example presented in \S \ref{ssec:numerical}.
Question (ii), instead, can be
addressed independently of (i) in the sense that, no matter how many basis functions have to be computed in order to meet
a given target accuracy, can one come up with methods that guarantee generating a {\em nearly minimal number} of such basis functions?
This has to do with {\em how to sample} the solution manifold and is the central theme in this paper.

The most prominent way of generating the reduced bases is a certain
{\em greedy sampling} of the manifold $\cM$.
Contriving {\em greedy sampling strategies} that give rise to reduced bases of
nearly minimal length, in a sense to be made precise below, also for {\em non-coercive or unsymmetric singularly perturbed problems} is the central objective in this paper.
We remark though that a greedy parameter search in its standard form is perhaps not suitable for very high dimensional
parameter spaces without taking additional structural features of the problem into account.
The subsequent discussions do therefore not target specifically the large amount of recent work on stochastic elliptic PDEs, since
while greedy concepts are in principle well understood for elliptic problems they are per se not necessarily adequate for
infinitely many parameters without exploiting specific problem dependent structural information.

First, we recall in \S \ref{sec:greedy} a {\em greedy space growth} paradigm commonly used in all established RBMs.
To measure its performance in the sense of (ii) we follow \cite{Buffaetal} and compare the corresponding distances ${\rm dist}_U(\cM,U_n)$ to the smallest possible distances achievable by linear spaces
of dimension $n$, called {\em Kolmogorov $n$-widths}.
The fact that for {\em elliptic problems} the convergence rates for the greedy errors are essentially those of the $n$-width, and hence {\em rate-optimal},
is shown in \S \ref{sect:growth} to be ultimately reduced to analyzing so called {\em weak greedy algorithms} in Hilbert spaces, see also \cite{BCDDPW,DPW-gr}.
However, for indefinite or strongly unsymmetric and singularly perturbed problems this
method usually operates far from optimality. We explain why this is the case and describe in \S \ref{sec:proj} a remedy proposed in \cite{DPW}.
A pivotal role is played by certain {\em well-conditioned variational formulations} of \eqref{1.1} which are then shown to
lead again to an optimal {\em outer greedy} sampling strategy also for non-elliptic problems. An essential additional
ingredient consists of certain stabilizing
{\em inner greedy loops}, see \S \ref{sec:d-greedy}. The obtained rate-optimal scheme is illustrated by a numerical example
addressing convection dominated convection-diffusion problems in
\S \ref{ssec:numerical}. We conclude in \S \ref{sec:duality} with applying these concepts to the efficient evaluation of
quantities of interest.

%%%%%%%%%%%%%%%%%%%%%%%%%%%%%
\section{The Greedy Paradigm}\label{sec:greedy}
%%%%%%%%%%%%%%%%%%%%%%%%%%%%%
The by far most prominent strategy for constructing reduced bases for a given parameter dependent problem \eqref{1.1} is
the following greedy procedure, see e.g. \cite{RHP}. The basic idea is that, having already constructed a reduced space $U_n$ of dimension
$n$, find an element $u_{n+1}=u(\mup_{n+1})$ in $\cM$ that is farthest away from the current space $U_n$, i.e.,
that maximizes the best approximation error from $U_n$ and then grow $U_n$ by setting $U_{n+1}:= U_n + \Span\,\{u_{n+1}\}$. Hence, denoting
by $P_{U,U_n}$ the $U$-orthogonal projection onto $U_n$,
\begin{equation}
\label{ideal}
\mup_{n+1}:= \argmax_{\mup\in \cY}\|u(\mup)- P_{U,U_n}u(\mup)\|_U, \quad u_{n+1}:= u(\mup_{n+1}).
\end{equation}
Unfortunately, determining such an exact maximizer is computationally way too expensive even in an offline phase because one would have to compute
for a sufficiently dense sampling of $\cY$ the exact solution $u(\mup)$ of \eqref{1.1} in $U$ (in practice in $U_\cN$). Instead one
tries to construct more efficiently computable {\em surrogates} $R(\mup,U_n)$ satisfying
\begin{equation}
\label{surrogate}
\|u(\mup)- P_{U,U_n}u(\mup)\|_U\le R(\mup,U_n), \quad \mup\in \cY.
\end{equation}
Recall that ``efficiently computable'' in the sense of offline-feasibility means that for each $\mup\in\cY$, the surrogate
$R(\mup,U_n)$ can be evaluated by solving only a problem of size $n$ in the reduced space $U_n$. Deferring an explanation of the
nature of such surrogates, Algorithm \ref{alg:greedy} described below is a
typical offline-feasible {\em surrogate based greedy algorithm} (\textsc{SGA}). Clearly, the maximizer in \eqref{greedy1} below is not necessarily unique. In case several maximizers exist it does not matter which one is selected.

\begin{algorithm}[htb]
  \caption{surrogate based greedy algorithm}
  \label{alg:greedy}
  \begin{algorithmic}[1]
    \Function{\textsc{SGA}}{}
    \State Set $U_0 := \{0\}$, $n=0$,
    \While{$\argmax_{\mup\in\cY} R(\mup,U_n) \ge tol$}
      \State
      \begin{equation}
      \begin{aligned}
      \mup_{n+1} & := \argmax_{\mup\in\cY} R(\mup,U_n), \\
      u_{n+1} & := u(\mup_{n+1}), \\
      U_{n+1}& := {\rm span}\,\big\{U_n,\{u(\mup_{n+1})\}\big\} = {\rm span}\,\{u_1, \dots, u_{n+1}\}
      \end{aligned}
      \label{greedy1}
      \end{equation}
    \EndWhile
    \EndFunction
  \end{algorithmic}
\end{algorithm}

Strictly speaking, the scheme {\SGA} is still idealized since:\\[2mm]
(a) computations cannot be carried out in $U$;\\[2mm]
(b) one cannot parse through all of a continuum $\cY$ to maximize $R(\mup,U_n)$.\\[0.5mm]

Concerning (a), as mentioned earlier computations in $U$ are to be understood as synonymous to computations
in a sufficiently large truth space $U_\cN$ satisfying all targeted accuracy tolerances for the underlying application.
Solving problems in $U_\cN$ is strictly confined to the offline phase and the number of such solves should remain of the order
of $n = {\rm dim}\,U_n$.
We will not distinguish in what follows between $U$ and $U_\cN$ unless such a distinction matters.

As for (b), the maximization is usually performed with the aid of a complete search over a {\em discrete} subset of $\cY$.
Again, we will not distinguish between a possibly continuous parameter set and a suitable training subset.
In fact, continuous optimization methods that would avoid a complete search have so far not proven to work well
since each greedy step increases the number of local maxima of the objective functional. Now, how fine such
a discretization for a complete search should be depends on how smoothly the $u(\mup)$ depend on $\mup$.
But even when such a dependence is very smooth a coarse discretization of a {\em high-dimensional} parameter set
$\cY$ would render a complete search infeasible so that, depending on the problem at hand,
one has to resort to alternate strategies such as, for instance, random sampling. However,
since it seems that (b) can only be answered for a specific problem class we will not address this issue in
this paper any further.

Instead, we focus on general principles which guarantee the following. Loosely speaking the reduced
spaces based on sampling $\cM$ should perform {\em optimally} in the sense that the resulting spaces $U_n$
have the (near) ``smallest dimension'' needed to satisfy a given target tolerance while the involved offline and online
cost remains feasible in the sense indicated above.
To explain first what is meant by ``optimal'' let us denote the {\em greedy error} produced by {\SGA} as
\begin{equation}
\label{sigman}
\sigma_n(\cM)_U := \max_{v\in \cM}\inf_{\bar u\in U_n}\|v-\bar u\|_U = \max_{\mup\in\cY}\|u(\mup)-P_{U,U_n} u(\mup)\|_U.
\end{equation}
Note that if we replaced in \eqref{sigman} the space $U_n$ by {\em any} linear subspace $W_n\subset U$ and infimize
the resulting distortion over {\em all} subspaces of $U$ of dimension at most $n$, one obtains the classical
{\em Kolmogorov $n$-widths} $d_n(\cM)_U$ quantifying the ``thickness'' of a compact set, see \eqref{Kolmo}. One trivially has
\begin{equation}
\label{lower}
d_n(\cM)_U \leq \sigma_n(\cM)_U,\quad n\in\N.
\end{equation}
Of course, it would be best if one could reverse the above inequality. We will discuss in the next section to what extent
this is possible.

To prepare for such a discussion we need more information about how the surrogate $R(\mup,U_n)$ relates to
the actual error $\|u(\mup)-P_{U,U_n}u(\mup)\|_U$ because the surrogate drives
the greedy search and one expects that the quality of the snapshots found in {\SGA} depends on
how ``tight'' the upper bound in \eqref{surrogate} is.

To identify next the essential conditions on a ``good'' surrogate it is instructive
to consider the case of {\em elliptic} problems. To this end,
suppose that
\begin{equation*}
\langle B_\mup u,v\rangle = b_\mup(u,v) =\langle f,v\rangle,\quad u,v\in U,
\end{equation*}
is a uniformly $U$-coercive bounded bilinear form and $f \in U'$, i.e., there exist constants $0< c_1\le C_1 <\infty$
such that
\begin{equation}
\label{ellip}
c_1 \|v\|_U^2\leq b_\mup(v,v), \quad | b_\mup(u,v)| \leq C_1\|u\|_U \|v\|_U,\,\, u,v \in U,\, \mup\in \cY,
\end{equation}
holds uniformly in $\mup\in \cY$. The operator equation \eqref{1.1} is then equivalent to:
given $f\in U'$ and a $\mup\in \cY$, find $u(\mup)\in U$ such that
\begin{equation}
\label{ell-prob}
b_\mup(u(\mup),v)=\langle f ,v\rangle, \quad v\in U.
\end{equation}
Ellipticity has two important well-known consequences. First, since \eqref{ellip} implies
$\|B_\mup\|_{U\to U'} \le C_1$, $\|B_\mup^{-1}\|_{U'\to U}\le c_1^{-1}$ the operator $B_\mup : U \to U'$ has a finite condition number
\begin{equation}
\label{condnumber}
\kappa_{U,U'}(B_\mup):= \|B_\mup\|_{U\to U'}\|B_\mup^{-1}\|_{U'\to U}\le C_1/c_1
\end{equation}
which, in particular, means that residuals in $U'$ are uniformly comparable to errors in $U$
\begin{equation}
\label{error-residual}
C_1^{-1}\|f - B_\mup \bar u\|_{U'}\leq \|u(\mup)- \bar u\|_U \le c_1^{-1} \|f - B_\mup \bar u\|_{U'}, \quad \bar u\in U,\, \mup\in \cY.
\end{equation}
Second, by C\'{e}a's Lemma, the Galerkin projection $\Pi_{\mup,U_n}$ onto $U_n$ is up to a constant as good as the {\em best approximation}, i.e.,
under the assumption \eqref{ellip}
\begin{equation}
\label{Gal}
\|u(\mup)- \Pi_{\mup,U_n}u(\mup)\|_U \le {\frac{C_1}{c_1}}\inf_{v\in U_n}\|u(\mup)-v\|_U.
\end{equation}
(When $b_\mup(\cdot,\cdot)$ is in addition symmetric $C_1/c_1$ can be replaced by $(C_1/c_1)^{1/2}$.)
Hence, by \eqref{error-residual} and \eqref{Gal},
\begin{equation}
\label{res-surr}
R(\mup,U_n):= c_1^{-1} \sup_{v\in U}\frac{\langle f, v\rangle -
b_\mup(\Pi_{\mup,U_n}u(\mup),v)}{\|v\|_U}
\end{equation}
statisfies more than just \eqref{surrogate}, namely it provides also a uniform {\em lower bound}
\begin{equation}
\label{surrogate-2}
\frac{c_1}{C_1}R(\mup,U_n)\le \|u(\mup)- P_{U,U_n}u(\mup)\|_U, \quad \mup\in \cY.
\end{equation}
Finally, suppose that $b_\mup(\cdot,\cdot)$ depends {\em affinely} on the parameters in the sense that
\begin{equation}
\label{affine}
b_\mup(u,v)= \sum_{k=1}^M \theta_k(\mup)b_k(u,v),
\end{equation}
where the $\theta_k$ are smooth functions of $\mup\in\cY$ and the bilinear forms $b_k(\cdot,\cdot)$
are independent of $\mup$. Then, based on
suitable precomputations (in $U_\cN$) in the offline phase, the computation of $\Pi_{\mup,U_n}u(\mup)$ reduces
for each $\mup\in\cY$ to the solution of a rapidly assembled $(n\times n)$-system,
and $R(\mup,U_n)$ can indeed be computed very efficiently, see \cite{RHP,Grepl1,SVHDNP}.

An essential consequence of \eqref{surrogate} and \eqref{surrogate-2} can be formulated as follows.

\begin{prop}
\label{prop:weakgreedy}
Given $U_n\subset U$, the function $u_{n+1}$ generated by \eqref{greedy1} for $R(\mup,U_n)$ defined by \eqref{res-surr},
has the property that
\begin{equation}
\label{weakgreedy}
  \|u_{n+1}-P_{U,U_n}u_{n+1}\|_U \geq \frac{c_1}{C_1} \max_{v\in \cM}\min_{\bar u\in U_n}\|v-\bar u\|_U .
\end{equation}
\end{prop}
Hence, maximizing the residual based surrogate $R(\mup,U_n)$ (over a suitable discretization of $\cY$) is a computationally feasible way
of determining, up to a fixed factor $\gamma := c_1/C_1 \le 1$, the maximal distance between $\cM$ and $U_n$ and performs in this sense
almost as well as the ``ideal'' but computationally infeasible surrogate $R^*(\mu,U_n):= \|u(\mup)- P_{U,U_n}u(\mup)\|_U$.\\

\noindent
{\bf Proof of Proposition \ref{prop:weakgreedy}:} Suppose that $\bar\mup = \argmax_{\mup\in \cY}R(\mup,U_n)$, $\mup^*:= \argmax_{\mup\in\cY}\|u(\mup)-P_{U,U_n}u(\mup)\|_U$ so that
$u_{n+1}= u(\bar\mup)$.
Then, keeping \eqref{surrogate-2} and \eqref{Gal} in mind, we have
\begin{eqnarray*}
 \| u_{n+1}-P_{U,U_n}u_{n+1}\|_U &=& \| u(\bar\mup) -P_{U,U_n}u(\bar\mup)\|_U
 \ge {\frac{c_1}{C_1}} R(\bar\mup,U_n) \ge \frac{c_1}{C_1} R( \mup^*,U_n) \\
&\ge & \frac{c_1}{C_1} \|u(\mup^*) - P_{U,U_n}u(\mup^*)\|_U = \frac{c_1}{C_1} \max_{\mup\in\cY}\|u(\mup)- P_{U,U_n}u(\mup)\|_U,
\end{eqnarray*}
where we have used \eqref{surrogate} in the second but last step.
This confirms the claim. \hfill $\Box$\\

Property \eqref{weakgreedy} turns out to play a key role in the analysis of the performance of the scheme {\SGA}.

%%%%%%%%%%%%%%%%%%%%%%%%%%%%%
\section{Greedy Space Growth}\label{sect:growth}
%%%%%%%%%%%%%%%%%%%%%%%%%%%%%

Proposition \ref{prop:weakgreedy} allows us to view the algorithm {\SGA} as a special instance of the following scenario.
Given a compact subset $\cK$ of a Hilbert space $H$ with inner product $(\cdot,\cdot)$ inducing the norm $\|\cdot\|^2=(\cdot,\cdot)$,
consider the {\em weak greedy} Algorithm \ref{alg:wgreedy} ({\WGA}) below.

\begin{algorithm}[htb]
  \caption{weak greedy algorithm}
  \label{alg:wgreedy}
  \begin{algorithmic}[1]
    \Function{\textsc{WGA}}{}
    \State Set $H_0 := \{0\}$, $n=0$, $u_0:=0$, fix any $0<\gamma \leq 1$,
    \State given $H_n$, choose some $u_{n+1}\in \cK$ for which
\begin{equation}
\label{weakgreedy2}
\min_{v_n\in H_n}\|v_n- u_{n+1}\| \ge \gamma \max_{v\in \cK}\min_{v_n\in U_n}\|v-v_n\| =: \gamma \sigma_n(\cK)_H,
\end{equation}
and set $H_{n+1}:= H_n + \Span\,\{u_{n+1}\}$.
  \EndFunction
  \end{algorithmic}
\end{algorithm}

Note that again the choice of $u_{n+1}$ is not necessarily unique and what follows holds for {\em any} choice satisfying \eqref{weakgreedy2}.

Greedy strategies have been used in numerous contexts and variants. The current version is not to be confused though
with the {\em weak orthogonal greedy algorithm} introduced in \cite{T} for {\em approximating} a {\em function} by a linear combination of
$n$ terms from a {\em given} dictionary. In contrast, the scheme \WGA ~described in Algorithm \ref{alg:wgreedy} aims at {\em constructing}
a (problem dependent) dictionary with the aid of a PDE model. While greedy function approximation is naturally compared with
the {\em best $n$-term approximation} from the underlying dictionary (see \cite{BCDD,T} for related results), a
natural question here is to compare the corresponding greedy errors
\begin{equation*}
\sigma_n(\cK)_H := \max_{v\in \cK}\min_{v_n\in U_n}\|v-v_n\|=: { \dist}_H(\cK,U_n)
\end{equation*}
incurred when approximating a compact set $\cK$
with the smallest possible deviation of $\cK$ from any $n$-dimensional linear space, given by the Kolmogorov $n$-widths
\begin{equation}
\label{Kolmo}
d_n(\cK)_H := \inf_{{\rm dim}V=n}\sup_{v\in \cK}\inf_{v_n\in V}\|v-v_n\| = \inf_{{\rm dim}V=n}{\dist}_H (\cK,V),
\end{equation}
mentioned earlier in the preceding section.
One trivially has $d_n(\cK)_H\le \sigma_n(\cK)_H$ for all $n\in\N$ and the question arises whether there actually exists a
constant $C$ such that
\begin{equation}
\label{compare1}
\sigma_n(\cK)_H\le C d_n(\cK)_H,\quad n\in \N.
\end{equation}
One may doubt such a relation
to hold for several reasons. First, orthogonal greedy {\em function approximation} performs in a way comparable to best $n$-term approximation
only under rather strong assumptions on the underlying given dictionary. Intutitively, one expects that errors made early on in the iteration are
generally hard to correct later although this intuition turns out to be misleading in the case of the present {\em set approximation}.
Second, the spaces $U_n$ generated by the greedy growth are restricted by being generated
only from snapshots in $\cK$ while the best spaces can be chosen freely, see the related discussion in \cite{BCDDPW}.

The comparison \eqref{compare1} was addressed first in \cite{Buffaetal} for the ideal case $\gamma =1$. In this case a bound of
the form $\sigma_n(\cK)_H\le Cn 2^n d_n(\cK)_H$ could be established for some absolute constant $C$.
This is useful only for cases where the $n$-widths decay faster than $n^{-1}2^{-n}$ which indeed turns out to
be possible for elliptic problems \eqref{ell-prob} with a sufficiently smooth affine parameter dependence \eqref{affine}.
In fact, in such a case the $u(\mup)$ can be even shown to be {\em analytic} as a function of $\mup$, see \cite{CDS1} and
the literature cited there.
It was then shown in \cite{BCDDPW} that the slightly better bound
\begin{equation}
\label{sharp}
\sigma_n(\cK)_H\le \frac{2^{n+1}}{\sqrt{3}} d_n(\cK)_H,\quad n\in \N,
\end{equation}
holds. More importantly, these bounds cannot be improved in general. Moreover, the possible exponential loss in accuracy
is not due to the fact the greedy spaces are generated by snapshots from $\cK$. In fact, denoting by $\bar d_n(\cK)_H$
the restricted ``inner'' widths, obtained by allowing only subspaces spanned by snapshots of $\cK$ in the competition, one can
prove that $\bar d_n(\cK)_H\le n d_n(\cK)_H$, $n\in\N$, which is also sharp in general \cite{BCDDPW}.

While these findings may be interpreted as limiting the use of reduced bases generated in a greedy fashion to problems
where the $n$-widths decay exponentially fast the situation turns out to be far less dim if one does {\em not} insist on a
{\em direct comparison} of the type \eqref{compare1} with $n$ being {\em the same} on both sides of the inequality. In \cite{BCDDPW,DPW-gr} the question is addressed whether a certain
{\em convergence rate} of the $n$-widths $d_n(\cK)_H$ implies some convergence rate of the greedy errors $\sigma_n(\cK)_H$.
The following result from \cite{BCDDPW} gave a first affirmative answer.

\begin{theorem}
\label{thm:poldecay} Let $0<\gamma\leq 1$ be the parameter in \eqref{weakgreedy2} and assume that
$d_0(\cK)_H\le M$ for some $M>0$. Then
\begin{equation*}
d_n(\cK)_H\le Mn^{-\alpha},\;\quad n\in \N,
\end{equation*}
for some $\alpha>0$, implies
\begin{equation}
\label{get}
\sigma_n(\cK)_H\le C Mn^{-\alpha}, \;\; n >0,
\end{equation}
where $C:= q^{\frac 1 2}(4q)^\alpha$ and
$q:= \lceil 2^{\alpha+1}\gamma^{-1}\rceil^2$.
\end{theorem}

This means that the weak greedy scheme may still be highly profitable even when the $n$-widths do not decay
exponentially. Moreover, as expected, the closer the weakness parameter $\gamma$ is to one, the better, which will
later guide the sampling strategies for constructing reduced bases.

Results of the above type are not confined to polynomial rates. A sub-exponential decay of the $d_n(\cK)_H$
with a rate $e^{-cn^\alpha}$, $\alpha \le1$
is shown in \cite{BCDDPW} to imply a rate
\begin{equation}
\label{subexp}
\sigma_n(\cK)_H \leq C(\alpha,\gamma)e^{-\tilde c n^{\tilde\alpha}},\quad \tilde\alpha =\alpha/(1+\alpha) ,\quad n\in\N.
\end{equation}
The principle behind the estimates \eqref{get}, \eqref{subexp} is to exploit a ``flatness'' effect or
what one may call ``conditional delayed comparison''. More precisely, given any $\theta \in (0,1)$ and defining
$q:= \lceil 2(\gamma\theta)\rceil^2$, one can show that (\cite[Lemma 2.2]{BCDDPW})
\begin{equation*}
\sigma_{n+qm}(\cK)_H \ge \theta \sigma_n(\cK)_H \quad \Rightarrow\quad \sigma_n(\cK)_H \leq q^{1/2} d_{m}(\cK)_H,\quad n\in \N.
\end{equation*}
Thus, a comparison between greedy errors and $n$-widths is possible when the greedy errors do not decay too quickly.
This is behind the diminished exponent $\tilde\alpha$ in \eqref{subexp}.

These results have been improved upon
in \cite{DPW-gr} in several ways employing different techniques yielding improved comparisons.
Abbreviating $\sigma_n:= \sigma_n(\cK)_H, d_n:= d_n(\cK)_H$,
a central result in the present general Hilbert space context states
that for any $N\ge 0, K\ge 1$, $1\le m<K$ one has
\begin{equation}
\label{general}
\prod_{i=1}^K \sigma_{N+i}^2 \le \gamma^{-2K}\Big(\frac{K}{M}\Big)^m\Big(\frac{K}{K-m}\Big)^{K-m}\sigma_{N+1}^{2m}d_m^{2(K-m)}.
\end{equation}
As a first important consequence, one derives from these inequalities a nearly direct comparison between $\sigma_n$ and $d_n$
without any constraint on the decay of $\sigma_n$ or $d_n$. In fact, taking $N=0, K=n$, and any $1\le m<n$ in \eqref{general},
using the monotonicity of the $\sigma_n$, one shows that
$\sigma_n^{2n}\le \gamma^{-2n}\Big(\frac nm\Big)^m\Big(\frac{n}{n-m}\Big)^{n-m} d_m^{2(n-m)}$ from which one deduces
\begin{equation}
\label{direct}
\sigma_n\le \sqrt{2}\gamma^{-1} \min_{1\le m <n}d_m^{\frac{n-m}{n}} ,\quad n\in \N.
\end{equation}
This, in particular, gives the direct unconditional comparison
\begin{equation*}
\sigma_{2n}(\cK)_H \leq \gamma^{-1}\sqrt{2 d_n(\cK)_H}, \quad n\in\N.
\end{equation*}
The estimate \eqref{direct} is then used in \cite{DPW-gr} to
improve on \eqref{subexp} establishing the bounds
\begin{equation}
\label{subexp2}
d_n(\cK)_H \le C_0e^{-c_0n^\alpha}, \quad \Rightarrow\quad \sigma_n(\cK)_H \leq \sqrt{2C_0}\gamma^{-1} e^{-c_1n^\alpha},\quad\,n\in\N,
\end{equation}
i.e., the exponent $\alpha$ is preserved by the rate for the greedy errors.
Moreover, one can recover \eqref{get} from \eqref{general} (with different constants).

Although not needed in the present context
the second group of results in \cite{DPW-gr} should be mentioned that concerns the extension of the weak greedy algorithm {\WGA}
to {\em Banach} spaces $X$ in place of the Hilbert space $H$. Remarkably, a direct comparison between $\sigma_n(\cK)_X$ and $d_n(\cK)_X$
similar to \eqref{general} is also established in \cite{DPW-gr}. The counterpart to \eqref{direct} reads $\sigma_{2n}\le 2\gamma^{-1}\sqrt{nd_n}$
i.e., one looses a factor $\sqrt{n}$ which is shown, however, to be necessary in general.

All the above results show that the smaller the weakness parameter $\gamma$ the stronger the derogation of the rate of the
greedy errors in comparison with the $n$-widths.

%%%%%%%%%%%%%%%%%%%%%%%%%%%%%
\section{What are the Right Projections?}\label{sec:proj}
%%%%%%%%%%%%%%%%%%%%%%%%%%%%%
As shown by \eqref{get} and \eqref{subexp2}, the weak greedy algorithm \WGA ~realizes optimal rates for essentially all
ranges of interest. A natural question is under which circumstances a surrogate based greedy algorithm \SGA~ is in this sense
also {\em rate-optimal}, namely ensures the validity of \eqref{get} and \eqref{subexp2}. Obviously, this is precisely the case
when new snapshots generated through maximzing the surrogate have the {\em weak greedy property} \eqref{weakgreedy2}.
Note that Proposition \ref{prop:weakgreedy} says that the {\em residual based surrogate} \eqref{res-surr} in the case of {\em coercive problems}
does ensure
the weak-greedy property so that \SGA ~ is indeed rate-optimal for {coercive problems}. Note also that the weakness
parameter $\gamma = c_1/C_1$ is in this case the larger the smaller the condition number of the operator $B_\mup$ is, see \eqref{condnumber}.
Obviously, the key is that the surrogate not only yields an upper bound for the best approximation error but also, up to a constant, a lower
bound \eqref{surrogate-2}, and the more tightly the best approximation eror is sandwiched by the surrogate the better the
performance of \SGA. Therefore, even if the problem is coercive for a very small $\gamma = c_1/C_1$, as is the case for
convection dominated {\em convection-diffusion problems}, in view of the dependence of the bounds in \eqref{get} and \eqref{subexp2}
on $\gamma^{-1}$, one expects that the performance of a greedy search based on
\eqref{res-surr} degrades significantly.

In summary, as long as algoritm {\SGA} employs a {\em tight surrogate} in the sense
that
\begin{equation}
\label{tight}
c_SR(\mup,U_n)\le \inf_{v\in U_n}\|u(\mup)-v\|_U \le R(\mup,U_n),\quad \mup\in \cY,
\end{equation}
holds for some constant $c_S>0$, independent of $\mup\in\cY$,
algorithm {\SGA} is {\em rate-optimal} in the sense of \eqref{get}, \eqref{subexp2}, i.e., it essentially realizes the $n$-width rates over all ranges of interest,
see \cite{DPW}. We refer to $c_S^{-1}:= \kappa_n(R)$ as the
{\em condition} of the surrogate $R(\cdot,U_n)$.
In the RBM community
the constant $c_S^{-1}$ is essentially the {\em stability factor} which is usually computed along with an approximate reduced
solution. Clearly, the bounds in \S \ref{sect:growth} also show that the quantitative performance of \SGA ~is expected to be the better the smaller the condition of the surrogate, i.e., the larger $c_S$.

As shown so far, coercive problems with a small condition number $\kappa_{U,U'}(B_\mup)$ represent an ideal setting for RBM
and standard Galerkin projection combined with the {\em symmetric surrogate} \eqref{res-surr}, based on measuring the residual
in the dual norm $\|\cdot\|_{U'}$ of the ``error-norm'' $\|\cdot\|_U$, identifies rate-optimal snapshots for a greedy space growth.
Of course, this marks a small segment of relevant problems. Formally, one can still apply these projections
and surrogates for any variational problem \eqref{ell-prob}
for which a residual can be computed. However, in general, for indefinite or unsymmetric singularly perturbed
problems, the tightness relation \eqref{tight} may no longer hold for surrogates of the form \eqref{res-surr} or, if it holds the condition
$\kappa_n(R)$ becomes prohibitively large. In this latter case, the upper bound of the best approximation error is too loose to direct
the search for proper snapshots. A simple example is the {\em convection-diffusion} problem: for $f\in (H_0^1(\Omega))'$ find
$\in H_0^1(\Omega)$, $\Omega \subset \R^d$, such that
\begin{equation}
\label{condiff}
\e(\nabla u,\nabla v) + ({\vec{b}}\cdot\nabla u,v)+ (cu,v)=: b_\mup(u,v)= \langle f,v\rangle,\quad v\in H^1_0(\Omega),
\end{equation}
where, for instance, $\mup =(\e,\vec{b})\in \cY:= [\e_0,1]\times S^{d-1}$, $S^{d-1}$ the $(d-1)$-sphere.
\begin{rem}
\label{rem:cond}
It is well nown that when $c - \frac 12 {\rm div} \vec{b}\ge 0$ the problem \eqref{condiff} has for any $f\in H^{-1}(\Omega):= (H_0^1(\Omega))'$
a unique solution. Thus, for $U:= H^1_0(\Omega)$
\eqref{ellip} is still valid but with $\kappa_{U,U'}(B_\mup)\sim \e^{-1}$ which becomes arbitrarily large for a correspondingly
small diffusion lower bound $\e_0$.
\end{rem}
The standard scheme {\SGA} indeed no longer performs nearly as well as in the well conditioned case.
The situation is even less clear when $\e =0$ (with modified boundary conditions)
where no ``natural'' variational formulation suggests itself (we refer to \cite{DPW} for a detailed discussion of these examples).
Moreover, for {\em indefinite problems} the Galerkin projection does
generally perform like the best approximation which also adversily affects tightness of the standard symmetric residual based surrogate
\eqref{res-surr}.

Hence, to retain rate-optimality of \SGA ~also for the above mentioned extended range of problems one has to
find a better surrogate than the one based on the symmetric residual bound in \eqref{res-surr}.
We indicate in the next section that such better surrogates can indeed be obtained at affordable
computational cost for a wide range of problems through combining {\em Petrov-Galerkin projections}
with appropriate {\em unsymmetric} residual bounds. The approach can be viewed as
{\em preconditioning} the continuous problem already on the infinite dimensional level.

%%%%%%%%%%%%%%%%%%%%%%%%%%%%%
\subsection{Modifying the Variational Formulation}\label{ssec:stabvar}
%%%%%%%%%%%%%%%%%%%%%%%%%%%%%
We consider now a wider class of (not necessarily coercive) variational problems
\begin{equation}
\label{b}
b(u,v)=\langle f,v\rangle,\quad v\in V,
\end{equation}
where we assume at this point only for each $f\in V'$ the existence of a unique solution $u\in U$, i.e., the operator $B: U\to V'$, induced by $b(\cdot,\cdot)$, is bijective. This is well known to be equivalent to the validity of
\begin{equation}
\label{infsup}
\left\{
\begin{array}{l}\displaystyle
\inf_{w\in W}\sup_{v\in V}\frac{b(w,v)}{\|w\|_U \|v\|_V}\ge \beta,
\quad \sup_{v\in V}\sup_{w\in U}\frac{b(w,v)}{\|w\|_U \|v\|_V}\le C_b,\\[5mm] \mbox{for}
\, v\in V\,
\exists\, w\in W,\, \mbox{such that } b(w,v)\neq 0,
\end{array}\right.
\end{equation}
for some constants $\beta, C_b$. However, one then faces two principal obstructions regarding an RBM based on the scheme {\SGA}:

(a)
first, as in the case of \eqref{condiff} for small diffusion, $\kappa_{U,V'}(B)\le C_b/\beta$
could be very large so that the corresponding error-residual relation
\begin{equation}
\label{err-res1}
\|u - v\|_U \le \beta^{-1}\|f- Bv\|_{V'},\quad v\in U,
\end{equation}
renders a corresponding residual based surrogate ill-conditioned.

(b) When $b(\cdot,\cdot)$ is not coercive,
the Galerkin projection does, in general, not perform as well as the best approximation.

The following approach has been used in \cite{DPW} to address both (a) and (b).
The underlying basic principle is not new, see \cite{BaMo84}, and variants of it have been used for different purposes in different contexts such as
least squares finite element methods \cite{MMRS} and, more recently, in connection with {\em discontinuous Petrov Galerkin methods}
\cite{DHSW,DJ1,DJ2}. In the context of RBMs the concepts of {\em natural norms} goes sort of half way by sticking in the end to Galerkin projections \cite{SVHDNP}. This marks an essential distinction from the approach in \cite{DPW} discussed later below.

The idea is to change the topology of one of the spaces so as to (ideally) make the corresponding induced operator an {\em isometry}, see
also \cite{DHSW}. Following \cite{DPW}, fixing for instance, $\|\cdot\|_V$, one can define
\begin{equation}
\label{hatU}
\|w\|_{\hat U}:= \sup_{v\in V}\frac{b(w,v)}{\|v\|_V} = \|B w\|_{V'}, \quad w\in U,
\end{equation}
which means that one has for $Bu=f$
\begin{equation}
\label{err-res2}
\|u-w\|_{\hat U}= \|f- Bw\|_{V'},\quad w\in U,
\end{equation}
a perfect error-residual relation. It also means that, replacing $\|\cdot\|_{U}$ in \eqref{infsup} by $\|\cdot\|_{\hat U}$, yields
the inf-sup constant $\hat\beta =1$. Alternatively, fixing $\|\cdot\|_U$, one may set
\begin{equation}
\label{hatV}
\|v\|_{\hat V}:= \sup_{w\in U}\frac{b(w,v)}{\|w\|_U} = \|B^* v\|_{U'}, \quad v\in V,
\end{equation}
to again arrive at an isometry $B: U\to \hat V'$, meaning
\begin{equation}
\label{err-res3}
\|u-w\|_{U}= \|f- Bw\|_{\hat V'},\quad w\in U.
\end{equation}
Whether the norm for $U$ or for $V$ is prescribed depends on the problem at hand and we refer to \cite{cdw,DPW,DHSW}
for examples of both type.

Next note that for any subspace $W\subset U$ one has
\begin{equation}
\label{bestappr}
u_W= \argmin_{w\in W}\|u- w\|_{\hat U} = \argmin_{w\in W}\|f- Bw\|_{V'},
\end{equation}
and analogously for the pair $(U,\hat V)$, i.e., the best approximation in the $\hat U$-norm
is a {\em minimum residual solution} in the $V'$-norm.

To use residuals in $V'$ as surrogates requires fixing a suitable discrete projection for a given trial space. In general, in particular when
$V\neq U$, the Galerkin projection is no longer appropriate since inf-sup stability of the {\em infinite dimensional} problem is no longer
inherited by an arbitrary pair of {\em finite dimensional} trial and test spaces. To see which type of projection
would be ideal, denote by $R_U: U'\to U$ the Riesz map defined for any
linear functional $\ell \in U'$ by
\begin{equation*}
\langle \ell, w\rangle = (R_U\ell,w)_U,\quad w\in U.
\end{equation*}
Then, by \eqref{hatU}, for any $w\in W\subset U$, taking $v:= R_VBw\in V$ one has
\begin{equation*}
b(w,v)= \langle Bw, v\rangle = \langle Bw, R_VBw\rangle = (Bw,Bw)_{V'} = (w,w)_{\hat U}.
\end{equation*}
Thus, in particular,
\begin{equation*}
b(u-u_h,R_V b w)= (u-u_h,w)_{\hat U},
\end{equation*}
i.e., given $W\subset U$, using $V(W):= R_V B (W)$ as a test space in the {\em Petrov-Galerkin} scheme
\begin{equation}
\label{PG1}
b(u_h, v)= \langle f,v\rangle,\quad v\in V(W):= R_V B (W),
\end{equation}
is equivalent to computing the $\hat U$-orthogonal projection of the exact solution $u$ of \eqref{b} and hence the best $\hat U$-approximation
to $u$. One readily sees that this also means
\begin{equation}
\label{infsup2}
\inf_{w\in W} \sup_{v\in V(W)}\frac{b(w,v)}{\|w\|_{\hat U}\|v\|_V}=1,
\end{equation}
i.e., we have a Petrov-Galerkin scheme for the pair of spaces $W, V(W)$ with perfect stability and the Petrov-Galerkin projection
is the best $\hat U$-projection. Unfortunately, this is not of much help yet, because computing the
{\em ideal test space} $V(W) = R_V B(W) =B^{-*}R_{\hat U}^{-1}(W)$ is not numerically feasible. Nevertheless, it provides a useful orientation
for finding good and practically realisable pairs of trial and test spaces, as explained next.

%%%%%%%%%%%%%%%%%%%%%%%%%%%%%
\subsection{A Saddle Point Formulation}\label{ssec:saddle}
%%%%%%%%%%%%%%%%%%%%%%%%%%%%%
We briefly recall now from \cite{DHSW,DPW} an approach to deriving from the preceding observations a practically
feasible numerical scheme which, in particular, fits into the context of RBMs.
Taking \eqref{bestappr} as point of departure we notice that the minimization of $\|f- Bw\|_{V'}$ over $W$ is a least squares
problem whose normal equations read: find $u_W\in W$ such that (with $R_{V'}=R_V^{-1}$)
\begin{equation}
\label{normal}
0= (f- Bu_W,Bw)_{V'} = \langle R_V^{-1} (f- B u_W), Bw\rangle,\quad w\in W.
\end{equation}
Introducing the auxiliary variable $r := R_V^{-1} (f- B u_W)$ which is equivalent to
\begin{equation}
\label{R}
\langle R_V r, v\rangle =(r,v)_V = \langle f- B u_w,v\rangle,\quad v\in V,
\end{equation}
the two relations \eqref{normal} and \eqref{R} can be rewritten in form of the {\em saddle point problem}
\begin{equation}
\label{saddlepoint}
\begin{array}{lcll}
(r,v)_V + b(u_W,v) &=& \langle f,v\rangle,& v\in V.\\
b(w,r) &=& 0, & w\in W.
\end{array}
\end{equation}
The corresponding inf-sup constant is still one (since the supremum is taken over all of $V$) and $(\cdot,\cdot)_V$ is
a scalar product so that \eqref{saddlepoint} has a unique solution $u_W$, see e.g. \cite{BF}. Taking for any $w\in W$ the test function $v=R_VBw\in V$
in the first line of \eqref{saddlepoint}, one obtains
\begin{equation*}
(r,v)_V= (r,R_VBw)_V = \langle r,Bw\rangle= b(w,r)=0,
\end{equation*}
by the second
line in \eqref{saddlepoint} so we see that $\langle f,R_VBw\rangle = b(u_W,R_VBw)$ holds for all $w\in W$ which means that $u_W$
solves the ideal Petrov-Galerkin problem \eqref{PG1}. Thus \eqref{saddlepoint}
is equivalent to the ideal Petrov Galerkin scheme \eqref{PG1}.

Of course, \eqref{saddlepoint} is still not realizable since the space $V$ is still the full infinite dimensional
space. One more step to arrive at a realizable scheme is the following: given the finite dimensional space $W$ find a finite dimensional space
$Z\subset V$ so that when $V$ in \eqref{saddlepoint} is replaced by $Z$, one obtains a {\em stable} finite dimensional saddle point problem
which is the same as saying that its inf-sup constant is safely bounded away from zero. Since $Z=V$ would yield perfect stability the choice
of $Z\subset V$ can be viewed as a {\em stabilization}. To quantify this we follow \cite{DPW} and say that for some $\delta \in (0,1)$,
$Z\subset V$ is $\delta$-{\em proximal} for
$W\subset U$ if $Z$ is suffciently close to the ideal test space $V(W)=R_V B(W)$ in the sense that
\begin{equation}
\label{deltaprox}
\|(I - P_{V,Z})R_VB w\|_V \le \delta \|R_VB w\|_V,\quad w\in W.
\end{equation}
The related main findings from \cite{DPW} can be summarized as follows.
\begin{theorem}
\label{thm:saddle}
(i) The pair $(u_{W,Z},r_{W,Z})\in W\times Z \subset U\times V$ solves the saddle point problem
\begin{equation}
\label{saddlepoint2}
\begin{array}{lcll}
(r_{W,Z},v)_V + b(u_{W,Z},v) &=& \langle f,v\rangle,& v\in Z,\\
b(w,u_{W,Z}) &=& 0, & w\in W,
\end{array}
\end{equation}
if and only if $u_{W,Z}$ solves the Petrov-Galerkin problem
\begin{equation}
\label{PG2}
b(u_{W,Z},v)= \langle f, v\rangle , \quad v\in P_{V,Z}(R_V B(W)).
\end{equation}
(ii) If $Z$ is $\delta$-proximal for $W$, \eqref{PG2} is solvable and one has
\begin{equation}
\label{bestappr2}
\begin{array}{rcl}
\|u-u_{W,Z}\|_{\hat U} &\le &\frac{1}{1-\delta}\inf_{w\in W}\|u_{W,Z}-w\|_{\hat U},\\[2mm]
\|u-u_{W,Z}\|_{\hat U}
+ \|r_{W,Z}\|_V &\le & \frac{2}{1-\delta}\inf_{w\in W}\|u_{W,Z}-w\|_{\hat U}.
\end{array}
\end{equation}
(iii) $Z$ is $\delta$-proximal for $W$ if and only if
\begin{equation}
\label{infsupdelta}
\inf_{w\in W}\sup_{v\in Z}\frac{b(w,v)}{\|w\|_{\hat U}\|v\|_V}\ge \sqrt{1-\delta^2}.
\end{equation}
\end{theorem}
Note that \eqref{saddlepoint2} involves ordinary bilinear forms and finite dimensional spaces $W,Z$ and (iii) says that
the $V$-projection of the ideal test space $R_VB(W)$ onto $Z$ is a good test space if and only if
$Z$ is $\delta$-proximal for $W$. Loosely speaking,
$Z$ is large enough to ``see'' a substantial part of the ideal test space $R_VB(W)$ under projection.
The perhaps most important messages to be taken home regarding the RBM context read as follows.

\begin{rem}
\label{rem:takehome}
(i) The Petrov-Galerkin scheme \eqref{PG2} is realized
through the saddlepoint problem \eqref{saddlepoint2} {\em without} explicitly computing the test space $P_{V,Z}(R_V B(W))$.\\[1.5mm]
(ii)
Moreover, given $W$, by compactness and \eqref{deltaprox}, one can in principle enlarge $Z$
so as to make $\delta$ as small as possible, a fact that will be explointed later. \\[1.5mm]
(iii) The solution component $u_{W,Z}$ is a near-best approximation to the exact solution $u$ in the
$\hat U$-norm.\\[1.5mm]
(iv) $r_{W,Z}$ can be viewed as a {\em lifted residual} which tends to zero in $V$ when $W$ grows
and can be used for a-posteriori error estimation, see \cite{DHSW}. In the Reduced Basis context this can be exploited for certifying the accuracy of the
truth solutions and for constructing computationally feasible surrogates for the construction of the reduced bases.
\end{rem}

%%%%%%%%%%%%%%%%%%%%%%%%%%%%%
\section{The Reduced Basis Construction}\label{sec:d-greedy}
%%%%%%%%%%%%%%%%%%%%%%%%%%%%%
We point out next how to use the preceding results for sampling the solution manifold $\cM$ of a given
{\em paramteric family of variational problems}: given $\mup\in \cY$, $f\in V_\mup'$, find $u(\mup)\in U_\mup$ such that
\begin{equation}
\label{parb}
b_\mup(u(\mup),v)=\langle f,v\rangle,\quad v\in V_\mup,
\end{equation}
in a way that the corresponding subspaces are rate-optimal. We will always assume that the dependence of the bilinear
form $b_\mup(\cdot,\cdot)$ on $\mup\in \cY$ is affine in the sense
of \eqref{affine}.

As indicated by the notation the spaces $U_\mup, V_\mup$ for which
the variational problems are well-posed in the sense that the induced operator $B_\mup:U_\mup\to V_\mup'$
is bijective, could depend on $\mup$ through $\mup$-dependent norms.
However, to be able to speak of a ``solution manifold'' $\cM$ as a compact subset of some ``reference Hilbert space'', the
norms $\|\cdot\|_{U_\mup}$ should be {\em uniformly} equivalent to some {\em reference} norm $\|\cdot\|_U$ which has to be taken into account
when formulating \eqref{parb}.
In fact, under this condition, as shown in \cite{DPW}, for well-posed variational formulations of pure transport problems
the dependence of the test spaces
$V_\mup$ on $\mup\in\cY$
is essential, in that
\begin{equation}
\label{capV}
V:= \bigcap_{\mup\in \cY}V_\mup,
\end{equation}
is a strict subset of each individual $V_\mup$. This complicates the construction of
a tight surrogate. We refer to \cite{DPW} for ways of dealing with this obstruction and confine the subsequent discussion
for simplicity to cases where the test norms $\|\cdot\|_{V_\mup}$ are also uniformly equivalent to a single reference norm $\|\cdot\|_V$,
see the example later below.

Under the above assumptions, the findings of the preceding section will be used next to contrive a well-conditioned tight surrogate
even for non-coercice or severely ill-conditioned variational problems which is then in general unsymmetric, i.e., $V_\mup\neq U_\mup$.
These surrogates will then be used in \SGA.
To obtain such a residual based well-conditioned surrogate
in the sense of \eqref{tight}, we first {\em renorm} the pairs of spaces $U_\mup$ or $V_\mup$ according to \eqref{hatU} or \eqref{hatV}.
In anticipation of the example below, for definiteness we concentrate on \eqref{hatU} and refer to \cite{DPW} for a discussion of \eqref{hatV}.
As indicated above, we assume further that the norms $\|\cdot\|_{\hat U_\mup}, \|\cdot\|_{V_\mup}$ are equivalent to
reference norms $\|\cdot\|_{\hat U}, \|\cdot\|_V$, respectively.

%%%%%%%%%%%%%%%%%%%%%%%%%%%%%
\subsection{The Strategy}
%%%%%%%%%%%%%%%%%%%%%%%%%%%%%
Suppose that we have
already constructed a pair of spaces $U_n\subset U_\mup ,V_n \subset V_\mup$, $\mup\in\cY$, such that for a given $\delta<1$
\begin{equation}
\label{infsupmup}
\inf_{w\in U_n}\sup_{v\in V_n}\frac{b_\mup(w,v)}{\|w\|_{\hat U_\mup}\|v\|_{V_\mup}}\ge \sqrt{1-\delta^2},\quad \mup \in \cY,
\end{equation}
i.e., $V_n\subset V$ is $\delta$-proximal for $U_n\subset U$.
Thus, by Theorem \ref{thm:saddle}, the parametric saddle point problem
\begin{equation}
\label{saddlepoint3}
\begin{array}{lcll}
(r_{n}(\mup),v)_{V_\mup} + b_\mup(u_{n}(\mup),v) &=& \langle f,v\rangle,& v\in V_n,\\
b(w,r_{n}(\mup)) &=& 0, & w\in W,
\end{array}
\end{equation}
has for each $\mup\in\cY$ a unique solution $(u_n(\mup),r_n(\mup))\in U_n \times V_n$.
By the choice of norms we know that
\begin{equation}
\label{know1}
\|u(\mup)- u_n(\mup)\|_{\hat U_\mup}= \| f- B_\mu u_n(\mup)\|_{V_\mup'}, \quad \mup \in \cY,
\end{equation}
i.e.,
\begin{equation}
\label{surr3}
R(\mup,U_n\times V_n):= \| f- B_\mu u_n(\mup)\|_{V_\mup'}, \quad \mup \in \cY
\end{equation}
suggests itself as a surrogate. There are some subtle issues about how to evaluate $R(\mup,U_n\times V_n)$ in
the dual $V_\cN'$ of a sufficiently large truth space $V_\cN\subset V_\mup$, $\mup\in \cY$, so as to faithfully reflect
errors in $\hat U_\mu$, not only in the truth space $U_\cN\subset U_\mup$ but in $\hat U$, and
how these quantities are actually related to the auxiliary variable $\|r_n(\mup)\|_{V_\mup}$ which is
computed anyway. As indicated before, these issues are aggrivated when the norms $\|\cdot\|_{V_\mup}$ are {\em not} all equivalent to a single
reference norm. We refer to a corresponding detailed discussion in \cite[\S 5.1]{DPW} and
continue working here for simplicity with the idealized version \eqref{surr3} and assume its offline feasibility.

Thus, we can evaluate the errors $\|u(\mup)- u_n(\mup)\|_{\hat U_\mup}$ and can determine a maximizing parameter
$\mup_{n+1}$ for which
\begin{equation}
\label{max}
\|u(\mup_{n+1})- u_n(\mup_{n+1})\|_{\hat U_\mup} = \max_{\mup\in \cY} \| f- B_\mu u_n(\mup)\|_{V_\mup'}.
\end{equation}
Now relation \eqref{bestappr2} in Theorem \ref{thm:saddle} tells us that for each $\mup\in\cY$
\begin{equation}
\label{bestappr2a}
\|u(\mup)- u_n(\mup)\|_{\hat U_\mup}\leq (1-\delta)^{-1}\inf_{w\in U_n}\|u(\mup)-w\|_{\hat U_\mup},
\end{equation}
i.e., $u_n(\mup)$ is a near best approximation to $u$ from $U_n$ which is, in fact, the nearer to the best approximation the smaller $\delta$.
By \eqref{know1} and \eqref{bestappr2a}, the surrogate \eqref{surr3} is indeed well-conditioned with condition number close to one for small
$\delta$.

A natural strategy is now to enlarge $U_n$ to $U_{n+1}:= U_n +\Span\,\{u(\mup_{n+1})\}$.
In fact, this complies with the {\em weak gready} step \eqref{weakgreedy2}
in \S \ref{sect:growth} with weakness parameter $\gamma =(1-\delta)$ as close to one as one wishes, when $\delta$ is chosen accordingly small,
provided that the pair of spaces $U_n, V_n$ satisfies \eqref{infsupmup}.
A repetition would therefore, in principle, be a realization of Algorithm \ref{alg:greedy}, \SGA, establishing rate-optimality of this RBM.
Obviously, the critical condition for such a procedure to work is to ensure at each stage the validity of the weak-greedy condition \eqref{weakgreedy2}
which in the present situation means that the companion space $V_n$ is at each stage $\delta$-proximal for
$U_n$. So far we have not explained yet how to grow $V_n$ along with $U_n$ so as to ensure $\delta$-proximality.
This is explained in the subsequent section.

\begin{rem}
\label{rem:crucial}
One should note that, due to the possible parameter dependence of the norms $\|\cdot\|_{\hat U_\mup}, \|\cdot\|_{V_\mup}$
on $\mup$, obtaining tight surrogates with the aid of an explicit Petrov-Galerkin formulation,
would be infeasible in an RBM context because one would have to recompute the
corresponding (parameter dependent) test basis for each parameter query which is not online-feasible.
It is therefore actually crucial to employ the saddle point formulation in the context of RBMs since this allows us to determine
a space $V_n$ of somewhat larger dimension than $U_n$ but stabilizes the saddle point problem {\rm for all $\mup$ simultaneously}.
\end{rem}

%%%%%%%%%%%%%%%%%%%%%%%%%%%%%
\subsection{A Greedy Stabilization}
%%%%%%%%%%%%%%%%%%%%%%%%%%%%%
A natural option is to enlarge $V_n$ by the second component $r_n(\mup_{n+1})$ of \eqref{saddlepoint3}.
Note though that the lifted resilduals $r_n$ tend to zero as $n\to \infty$. Hence, the solution manifold
of the ($\mup$-dependent version of the) saddle point formulation \eqref{saddlepoint} has the form
\begin{equation*}
\cM \times \{0\},
\end{equation*}
where $\cM$ is the solution manifold of \eqref{parb} (since $r(\mup)=0$ for $\mup\in \cY$). Thus,
the spaces $V_n$ are {\em not} needed to approximate the solution manifold. Instead the sole purpose of the space $V_n$
is to guarantee stability. At any rate, the grown pair $U_{n+1}, V_n + \Span\,\{r_n(\mup_{n+1})\}=: V_{n+1}^0$ may fail to satisfy now
\eqref{infsupmup}.

Therefore, in general one has to further enrich $V_{n+1}^0$ by additional {\em stabilizing} elements again in a greedy fashion
until \eqref{infsupmup} holds for the desired $\delta$.
For problems that initially arise as natural saddle point problems such as the Stokes system, enrichments by so called {\em supremizers}
(to be defined in a moment)
have been proposed already in \cite{GV,GV2,Rozza}. In these cases it is possible to enrich $V_{n+1}^0$ by a {\em fixed}
a priori known number of such supremizers to guarantee inf-sup stability. As shown in \cite{DPW}, this is generally possible when
using fixed (parameter independent) reference norms $\|\cdot\|_{\hat U}$, $\|\cdot\|_V$ for $U$ and $V$. For the above more general scope of problems
a greedy strategy was proposed and analyzed in \cite{DPW}, a special case of which is also considered in \cite{GV2} without
analysis. The strategy in \cite{DPW} adds only as many
stabilizing elements as are actually needed to ensure stability and works for a much wider range of problems including
singularly perturbed ones. In cases where not all parameter dependent norms $\|\cdot\|_{V_\mup}$
are equivalent such a strategy is actually necessary and its convergence analysis is then more involved, see \cite{DPW}.

To explain the procedure, suppose that after growing $U_n$ to $U_{n+1}$ we have already generated an enrichment $V_{n+1}^k$
of $V_{n+1}^0$ (which could be, for instance, either $V_{n+1}^0:=V_n + \Span\,\{r_n(\mup_{n+1})\}$ or
$V_{n+1}^0:= V_n$) but the pair $U_{n+1}, V_{n+1}^k$ still fails to satisfy \eqref{infsupmup} for the given $\delta <1$.
To describe the next enrichment from $V_{n+1}^k$ to $V_{n+1}^{k+1}$
we first search for a parameter $\bar{\mup} \in \cY$ and a function $\bar{w} \in U_n$ for which the inf-sup condition \eqref{infsupmup} is worst, i.e.,
\begin{equation}
\label{arginfsup}
\sup_{v\in V_{n+1}^k}\frac{b_{\bar{\mup}}(\bar{w}, v)}{\|v\|_{V_{\bar{\mup}}}\|\bar{w}\|_{\hat U_{\bar{\mup}}}} = \inf_{\mup\in \cY} \left( \inf_{w\in U_{n+1}} \sup_{v\in V_{n+1}^k}\frac{b_\mup(w, v)}{\|v\|_{V_\mup}\|w\|_{\hat U_\mup}} \right).
\end{equation}
If this worst case inf-sup constant does not exceed yet $\sqrt{1-\delta^2}$, the current space $V_{n+1}^k$ does not contain
an effective supremizer {for $\bar{\mup}, \bar{w}$,} yet. However, since the truth space satisfies the uniform inf-sup condition \eqref{infsupmup} there {\em must exist} a good supremizer in the truth space which can be seen to be given by
\begin{equation}
\label{supremizer}
\bar{v} =
\argmax_{v \in V_{\bar{\mup}}}\frac{b_{\bar{\mup}}(\bar{w}, v)}{\|v\|_{V_{\bar{\mup}}} \|\bar{w}\|_{\hat U_{\bar{\mup}}}},
\end{equation}
providing the next enrichment
\begin{equation}
  V_{n+1}^{k+1} := {\rm span} \{V_{n+1}^k, \bar v\}.
  \label{eq:update}
\end{equation}
We defer some comments on the numerical realization of finding $\bar\mup, \bar v$ in \eqref{arginfsup}, \eqref{supremizer} to the next section.

This strategy can now be applied recursively until one reaches a satisfactory uniform inf-sup condition for the reduced spaces.
Again, the termination of this stabilization loop is easily ensured when \eqref{affine} holds and the norms $\|\cdot\|_{\hat U_\mup}$, $\|\cdot\|_{V_\mup}$
are uniformly equivalent to reference norms $\|\cdot\|_{\hat U }$, $\|\cdot\|_{V}$, respectively, but is more involved in the general case
\cite{DPW}.

%%%%%%%%%%%%%%%%%%%%%%%%%%%%%
\subsection{The Double Greedy Scheme and Main Result}
%%%%%%%%%%%%%%%%%%%%%%%%%%%%%
Thus, in summary, to ensure that the greedy scheme \SGA~ with the particular surrogate \eqref{surr3}, based on the corresponding {\em outer greedy} step for extending $U_n$ to $U_{n+1}$, has the {\em weak greedy property} \eqref{weakgreedy2}, one can employ
an {\em inner stabilizing greedy} loop producing a space $V_{n+1}=V_{n+1}^{k^*}$ which
is $\delta$-proximal for $U_{n+1}$. Here $k^*=k^*(\delta)$ is the number of enrichment steps needed to
guarantee the validity of \eqref{infsupmup} for the given $\delta$.
A sketchy version of the corresponding ``enriched'' \SGA, developed in \cite{DPW}, looks as follows:\\

\begin{algorithm}[htb]
  \caption{double greedy algorithm}
  \label{alg:doublegreedy}
  \begin{algorithmic}[1]
    \Function{\textsc{SGA-dou}}{}
  \State Initialize $U_1,V_1^0$, $\delta \in (0,1)$, target accuracy ${\rm tol}$, $n\leftarrow 1$,
  \While{$\sigma_n(\cM)> {\rm tol}$}
     \While{$U_n,V_n^0$ fail to satisfy \eqref{infsupmup}}
     \State compute $V_{n}$ with the aid of the inner stabilizing greedy loop,
     \EndWhile
    \State given $U_n,V_n$, satisfying \eqref{infsupmup},
    compute $U_{n+1}, V_{n+1}^0$ with the aid of the outer greedy \hspace*{10.5mm}step 4, \eqref{greedy1}
in algorithm {\SGA} for the surrogate \eqref{surr3},
       \EndWhile
  \EndFunction
  \end{algorithmic}
\end{algorithm}

%%%%%%%%%%%

As indicated above, both Algorithm \ref{alg:greedy}, \SGA~ and Algorithm \ref{alg:doublegreedy}, \DoubGA~ are
surrogate based greedy algorithms. The essential difference is that for non-coercive problems or problems with an originally  large
variational condition number in \DoubGA ~ an additional interior
greedy loop provides a tight well-conditioned (unsymmetric) surrogate which guarantees the desired
weak greedy property (with weakness constant $\gamma$ as close to one as one wishes) needed for rate-optimality.

Of course, the viability of Algorithm {\DoubGA} hinges mainly on two questions:

(a) how to find the worst inf-sup constant in \eqref{arginfsup} and how to compute the supremizer in \eqref{supremizer}?

(b) does the inner greedy loop terminate (early enough)?

As for (a), it is well known that, fixing bases for $U_{n},V_n^k$, finding the worst inf-sup constant amounts to determine
for $\mup\in\cY$ the cross-Gramian with respect to $b_\mup(\cdot,\cdot)$ and compute its smallest singular value. Since
these matrices are of size $n\times (n+k)$ and hence (presumably) of ``small'' size, a search over $\cY$ requires solving
only problems in the reduced spaces and are under the assumption
\eqref{affine} therefore offline-feasible. The determination of the corresponding supremizer $\bar v$ in \eqref{supremizer}, in turn,
is based on the well-known observation that
\begin{equation*}
 \argmax_{v \in V_{\bar{\mup}}}\frac{b_{\bar{\mup}}(\bar{w}, v)}{\|v\|_{V_{\bar{\mup}}} }
 = R_{V_{\bar{\mup}}} \bo_{\bar{\mup}} \bar{w},
\end{equation*}
which is equivalent to solving the Galerkin problem
\begin{equation*}
(\bar v, z)_{V_{\bar{\mup}}} = b_{\bar{\mup}}(\bar{w}, z),\quad z\in V_{\bar\mup}  .
\end{equation*}
Thus, each enrichment step requires one offline-Galerkin solve in the truth space.

A quantitative answer to question (b) is more involved. We are content here with a few related remarks and we refer to a detailed discussion of this issue in \cite{DPW}. As mentioned before, when all the norms $\|\cdot\|_{\hat U_\mup}, \|\cdot\|_{V_\mup}$, $\mup\in\cY$, are equivalent
to reference norms $\|\cdot\|_{\hat U}, \|\cdot\|_V$, respectively, the inner loop terminates after at most the number of terms in \eqref{affine}.
When the norms $\|\cdot\|_{V_\mup}$ are no longer uniformly equivalent to a single reference norm termination is less clear.
Of course, since all computations are done in a truth space which is finite dimensional, compactness guarantees termination after finitely many steps. However, the issue is that the number of steps should not depend on the truth space dimension. The reasoning in
\cite{DPW} used to show that (under mild assumptions) termination happens after a finite number of steps which is {\em independent}
of the truth space dimension, is based on the following fact. Defining $U^1_n(\mup):= \{w\in U_n: \|w\|_{\hat U_\mup}=1\}$, solving the problem
\begin{equation}
\label{delta-prox1}
(\bar \mup,\bar w):= \argmax_{\mup\in\cY;w\in U^1_n(\mup)} \inf_{\phi\in V_n^k}\|R_{V_\mup}B_\mup w-\phi\|_{V_\mup},
\end{equation}
when all the $\|\cdot\|_{\hat U_\mup}$-norms are equivalent to a single reference norm, can be shown to be equivalent to a greedy step of
the type \eqref{arginfsup} and can hence again be reduced
to similar small eigenvalue problems in the reduced space. Note, however, that \eqref{delta-prox1} is similar to a greedy space growth
used in the outer greedy loop and for which some understanding of convergence is available. Therefore, successive enrichments
based on \eqref{delta-prox1} are studied in \cite{DPW} regarding their convergence. The connection with the inner stabilizing loop based
on \eqref{arginfsup}
is that
\begin{equation*}
  \argmax_{\mup\in\cY;w\in U^1_n(\mup)} \inf_{\phi\in V_n^k}\|R_{V_{\bar \mup}}B_{\bar \mup} \bar w-\phi\|_{V_{\bar \mup}}\le \delta,
\end{equation*}
just means
\begin{equation*}
\inf_{\phi\in V_n^k}\|R_{V_{ \mup}}B_\mup w-\phi\|_{V_{ \mup}}\le \delta \|R_{V_{ \mup}}B_\mup\|_{V_\mup}=\delta \|w\|_{\hat U_\mup},\quad w\in U_n,\,\mup\in \cY,
\end{equation*}
which is a statement on $\delta$-proximality known to be equivalent to inf-sup stability, see Theorem \ref{thm:saddle},
and \eqref{deltaprox}.

A central result from \cite{DPW} can be formulated as follows, see \cite[Theorem 5.5]{DPW}.
\begin{theorem}
\label{thm:main}
If \eqref{affine} holds and the norms $\|\cdot\|_{\hat U_\mup}, \|\cdot\|_{V_\mup}$ are all equivalent to a single reference norm $\|\cdot\|_{\hat U},
\|\cdot\|_{V}$, respectively, and the surrogates \eqref{surr3}
are used, then the scheme {\DoubGA} is rate optimal, i.e., the greedy errors $\sigma_n(\cM)_{\hat U}$
decay at the same rate as the $n$-widths $d_n(\cM)_{\hat U}$, $n\to \infty$.
\end{theorem}

Recall that the quantitative behavior of the greedy error rates are directly related to those of the $n$-widths by $\gamma^{-1}= c_S^{-1}$,
see Theorem \ref{thm:poldecay}.
This suggests that a fast decay of $d_n(\cM)_{\hat U}$ is reflected by the corresponding greedy errors already for moderate
values of $n$ which is in the very interest of reduced order modeling. This will be confirmed by the
expamples below. In this context an important feature of {\DoubGA} is that through the choice of the $\delta$-proximility parameter
the weakness parameter $\gamma$ can be driven towards one, of course, at the expense of somewhat larger spaces $V_{n+1}$.
Hence, stability constants close to one are built into the method. This is to be contrasted by the conventional use of {\SGA}
based on surrogates that are not ensured to be well conditioned and for which the computation of the certifying stability constants
tends to be computationally expensive.

%%%%%%%%%%%%%%%%%%%%%%%%%%%%%
\subsection{A Numerical Example}\label{ssec:numerical}
%%%%%%%%%%%%%%%%%%%%%%%%%%%%%
The preceding theoretical results are illustrated next by a numerical example that brings out some of the main features of the scheme.
While the double greedy scheme applies to non-coercive or indefinite problems (e.g. see \cite{DPW} for pure transport)
we focus here on a classical {\em singularly perturbed} problem because it addresses also some principal issues for
RBMs regarding problems with {\em small scales}.
Specifically, we consider the {\em convection-diffusion} problem \eqref{condiff} on $\Omega = (0,1)^2$ for a simple {\em parameter
dependent convection} field
\begin{equation*}
\vec{b}(\mup):= \begin{pmatrix} \cos { \mup} \\ \sin { \mup} \end{pmatrix} ,\quad \mup \in [0,2\pi),\quad c=1,
\end{equation*}
keeping for simplicity the diffusion level $\e$ fixed but allowing it to be arbitrarily small.
All considerations apply as well to variable and parameter dependent diffusion with any arbitrarily small but strictly positive lower bound.
The ``transition''
to a pure transport problem is discussed in detail in \cite{DPW,W}.
Parameter dependent convection directions mark actually the more difficult case
and are, for instance, of interest with regard to kinetic models.

Let us first briefly recall the main challenges posed by \eqref{condiff}
for very small diffusion $\e$. The problem becomes obviously dominantly unsymmetric and singularly perturbed.
Recall that for each positive $\e$ the problem possesses for each $\mup\in \cY$ a unique solution $u(\mup)$ in $U=H^1_0(\Omega)$
that has a zero trace on the boundary $\partial\Omega$. However, as indicated earlier, the condition number $\kappa_{U,U'}(B_\mup)$
of the underlying convection-diffusion operator $B_\mup$, viewed as an operator from $U=H^1_0(\Omega)$ onto $U'=H^{-1}(\Omega)$,
behaves like $\e^{-1}$, that is, it becomes increasingly {\em ill conditioned}.
This has well known consequences for the performance of numerical solvers but above all for the stability of corresponding
discretizations.

We emphasize that the conventional mesh dependent stabilizations like SUPG (cf. \cite{hughesetal})
do {\em not} offer a definitive remedy because
the corresponding condition, although improved, remains very large for very small $\e$.
In \cite{PRcondiff} SUPG-stabilization for the offline truth calculations as well as for the low-dimensional online
Galerkin projections are discussed for moderate P\'{e}clet-numbers
of the order of up to $10^3$. In particular, comparisons are presented when only the offline phase uses stabilization while
the un-stabilized bilinear form is used in the online phase, see also the references in \cite{PRcondiff} for further related
work.

As indicated earlier, we also remark in passing that the singularly perturbed nature of the problem poses an additional difficulty concerning the
choice of the truth space $U_\cN$. In fact, when $\e$ becomes very small one may not be able to afford resolving correspondingly
thin layers in the truth space which increases the difficulty of capturing essential features of the solution by the reduced model.

This problem is addressed in \cite{DPW} by resorting to a weak formulation that does not use $H^1_0(\Omega)$ (or a renormed version
of it) as a trial space but builds on the results from \cite{cdw}.
A central idea is to enforce the boundary conditions on the outflow boundary $\Gamma_+(\mup)$ only weakly.
Here $\Gamma_+(\mup)$ is that portion of $\partial\Omega$ for which the inner product of the outward normal and
the convection direction is positive. Thus,
solutions are initially sought in the larger space $H^1_{0,\Gamma^-(\mup)}(\Omega)=: U_-(\mup)$ enforcing homogeneous boundary conditions
only on the {\em inflow} boundary $\Gamma_-(\mup)$. Since the outflow boundary, and hence also the inflow boundary depend on the parameter $\mup$,
this requires subdividing the parameter set into smaller sectors, here four, for which the outflow boundary $\Gamma_+=\Gamma_+(\mup)$ remains unchanged.
We refer in what follows for simplicity to one such sector denoted again by $\cY$.

\begin{figure}
  \centering
  \includegraphics[width=0.28\textwidth]{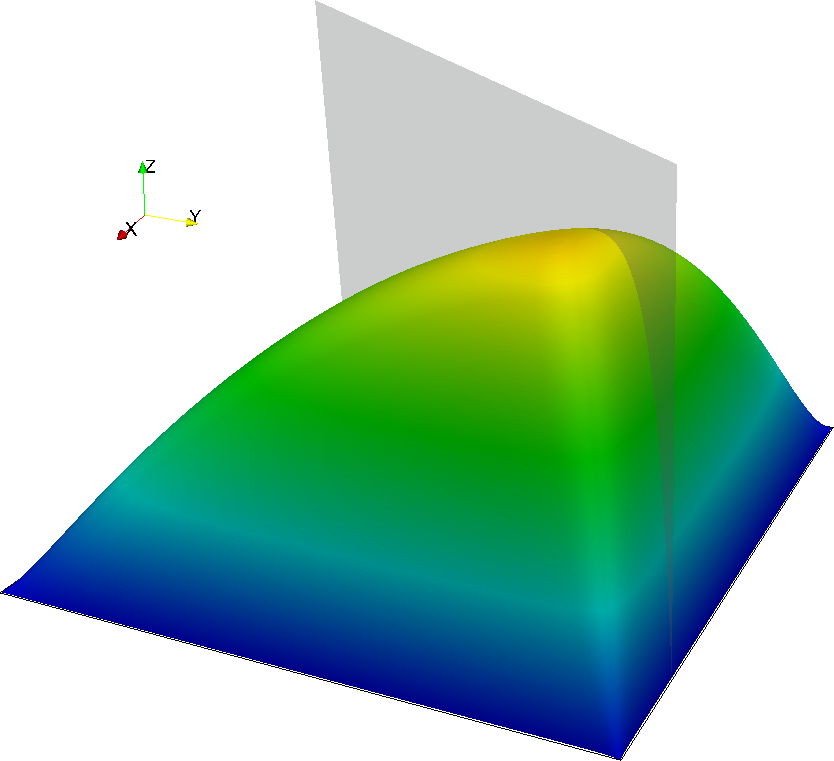}
  \qquad
  \hspace{1mm}
  \includegraphics[width=0.28\textwidth]{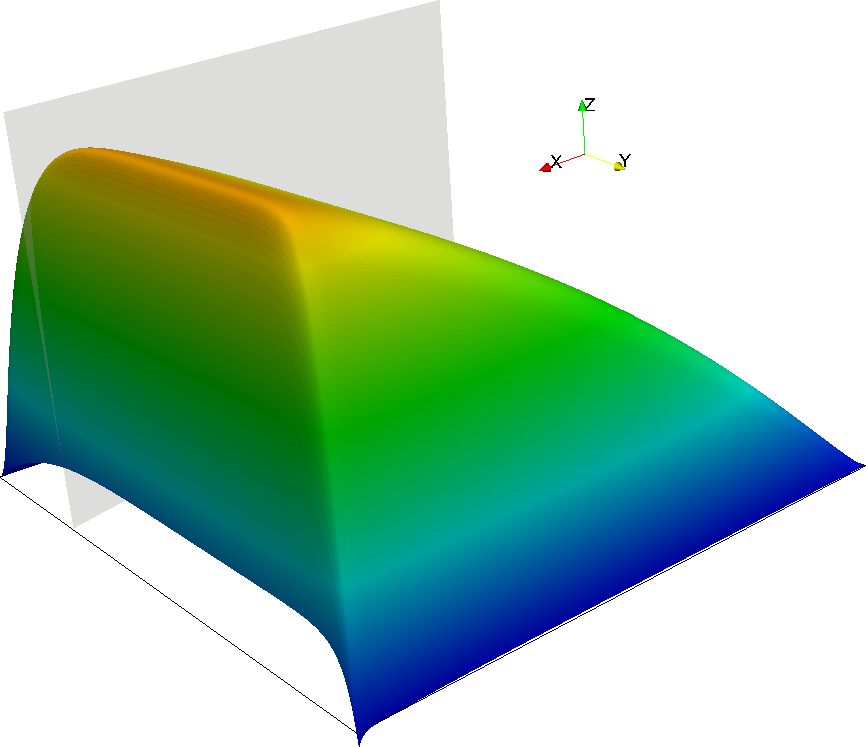}
  \qquad
  \includegraphics[width=0.28\textwidth]{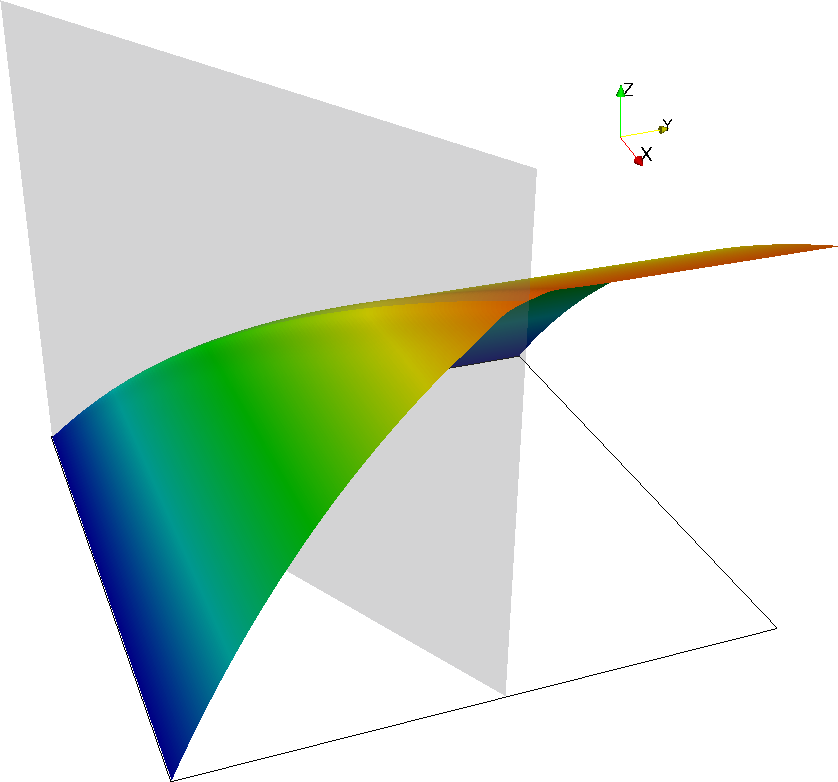}
  \caption{left: $\e= 2^{-5}, n= 6, n_V=13$; ~ middle:  $\e=2^{-7}, n=7, n_V=20$; ~ right: $\e= 2^{-26}, n=20,n_V=57$.}
  \label{figs:condiff}
\end{figure}

The following prescription of the {\em test space} falls into the category \eqref{hatU} where the norm for $U$ is adapted. Specifically,
choosing
\begin{eqnarray*}
s_\mup(u,v) &:= &\frac{1}{2}\big(\langle{B_\mup u, v}\rangle + \langle{B_\mup v, u}\rangle\big),\\
\|v\|_{V_\mup}^2& := & s_\mup(v, v) = \epsilon |v|_{H^1(\Omega)}^2 + \Big\|\Big(c - \frac{1}{2} {\rm div}\, \vec{b}(\mup)\Big)^{1/2} v \Big\|_{L_2(\Omega)}^2,
\end{eqnarray*}
in combination with a boundary penalization on $\Gamma_+$, we follow \cite{cdw,W} and define
\begin{equation*}
 \|u\|_{\bar{U}_\mup}^2 := \|\bar B_\mup u\|_{\bar V_\mup'}^2 = \| \bar B_\mup u\|_{V_\mup'}^2 + { \lambda} \|u\|_{H_b(\mu)}^2,
\end{equation*}
where $H_b(\mup)= H^{1/2}_{00}(\Gamma_+(\mup))$, $\bar V_\mup := V_\mup \times H_b(\mup)'$ and $\bar B_\mup$ denotes the
operator induced by this weak formulation over $\bar U_\mup := H^1_{0,\Gamma_-(\mup)}(\Omega)\times H_b(\mup)$.
The corresponding variational formulation is of minimum residual type (cf. \eqref{bestappr}) and reads
\begin{equation}
\label{penalize}
u(\mup)= {\rm argmin}_{w\in U_-(\mup)}\big\{\|\bar B_\mup w- f\|^2_{\bar V_\mup'} + \lambda\| w\|^2_{H_b(\mup)}\big\}.
\end{equation}
One can show that its (infinite dimensional) solution, whenever being sufficiently regular, solves also
the strong form of the convection diffusion problem \eqref{condiff}.
Figure \ref{figs:condiff} illustrates the effect of this formulation where we set $n= {\rm dim}\,U_n, n_V:= {\rm dim}\,V_n$

%%%%%%%%%%%%%%

\begin{figure}[ht]
  \begin{center}
  	\begin{tikzpicture}
      \begin{axis}[
          title={},
          ymin=0,
          xlabel={reduced basis trial dimension},
          width=0.5\textwidth
        ]
        \addplot table[x=dim_trial,y=max_residual]{all-cases-5};
      \end{axis}
    \end{tikzpicture}~
  \end{center}
  \caption{Convection-diffusion equation, $\varepsilon = 2^{-26}$, maximal a-posteriori error $0.00208994$}
  \label{fig:all-cases-5}
\end{figure}
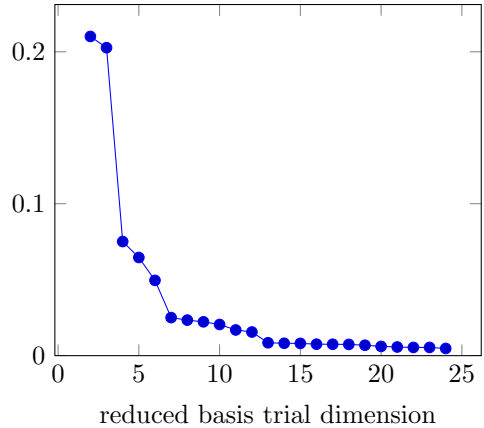

%%%%%%%%%%%%%%%%%

The shaded planes shown in Figure \ref{figs:condiff} indicate the convection direction for which the
snapshot is taken. For moderately large diffusion the boundary layer at $\Gamma_+$ is resolved by the truth space discretization and
the boundary conditions at the outflow boundary are satisfied exactly. For smaller diffusion in the middle example
the truth space discretization can no longer resolve the boundary layer and for very small diffusion (right) the solution is
close to the one for pure transport. The rationale of \eqref{penalize} is that all norms commonly used for convection diffusion equations
resemble the one chosen here, for instance in the form of a mesh dependent ``broken norm'', which means that most part of the
incurred error of an approximation is concentrated in the layer region, see e.g. \cite{Sangalli05,Verfurth05}. Hence, when the layers are not resolved by
the discretization, enforcing the boundary conditions does not improve accuracy and, on the contrary, may degrade
accuracy away from the layer by causing oscillations. The present formulation instead avoids any non-physical oscillations
and enhances accuracy in those parts of the domain where this is possible for the afforded discretization,
see \cite{cdw,DPW,W} for a detailed discussion. The following table quantifies the results for the case of small diffusion $\e =2^{-26}$
and a truth discretization whose a posteriori error bound is $0.002$.\\[3mm]

{\small
\begin{table}[h!]
  \caption{Convection-diffusion equation, $\varepsilon = 2^{-26}$, maximal a-posteriori error $0.00208994$}
  \begin{center}
    \begin{tabular}{|c|c|c|c|c||c|c|c|c|c|}
      \hline
      $n$ & $n_V$ & $\delta$ & surrogate & surr/a-post &$n$& $n_V$ & $\delta$ & surrogate & surr/a-post \\
      \hline \hline
      2 & 5 & 1.36e-03 & 2.10e-01 & 1.01e+02 &   14 & 39 & 1.17e-04 & 8.15e-03 & 3.90e+00 \\
      4 & 9 & 1.10e-02 & 7.51e-02 & 3.59e+01 & 16 & 45 & 9.79e-05 & 7.56e-03 & 3.62e+00 \\
      6 & 15 & 1.75e-03 & 4.95e-02 & 2.37e+01 & 18 & 51 & 6.32e-05 & 7.40e-03 & 3.54e+00 \\
      8 & 21 & 9.16e-04 & 2.34e-02 & 1.12e+01 & 20 & 57 & 4.74e-05 & 6.09e-03 & 2.92e+00 \\
      10 & 27 & 3.65e-04 & 2.05e-02 & 9.82e+00 &22 & 63 & 2.36e-05 & 5.43e-03 & 2.60e+00 \\
      12 & 33 & 3.34e-04 & 1.56e-02 & 7.45e+00 & 24 & 65 & 2.36e-05 & 4.73e-03 & 2.27e+00\\
      \hline
    \end{tabular}
  \end{center}
  \label{tab:all-cases-5}
\end{table}
}
\vspace*{2mm}

The columns $3$ and $8$ show the $\delta$ governing the condition of the saddle point problems (and hence
of the corresponding Petrov-Galerkin problems), see \eqref{infsupmup}, the greedy space growth is based upon.
Hence the surrogates are very tight giving rise to weakness parameters very close to one. As indicated in Remark \ref{rem:takehome}
one can use also an a posteriori bound for the truth solution based on the corresponding lifted residual. Columns $5$ and $10$
show therefore the relative accuracy of the current reduced model and the truth model.
This corresponds to the stability constants computed by conventional RBMs. Even for elliptic problems these latter ones are
significantly larger than the ones for the present singularly perturbed problem which are guaranteed to be close to one
by the method itself. Based on the a posteriori bounds for the truth solution (which are also obtained with the aid of tailored well-conditioned
variational formulations, see \cite{cdw}), the greedy space growth is stopped when
the surrogates reach the order of the truth accuracy. As illustrated in Figure \ref{fig:all-cases-5}, in the present example this is essentially
already the case for $\le 20$ trial reduced basis functions
and almost three times as many test functions. To show this ``saturation effect'' we have continued the space growth formally up to $n=24$
showing no further significant improvement which is in agreement with the resolution provided by the truth space. These relations agree with the theoretical predictions in \cite{DPW}. Figure \ref{fig:all-cases-5}
illustrates also the rapid gain of accuracy by the first few reduced basis functions which supports the fact that the solution manifold
is ``well seen'' by the Petrov-Galerkin surrogates.
More extensive numerical tests shown in \cite{DPW} show that the achieved stability is independent of the diffusion but the larger the diffusion
the smaller become the dimensions $n= {\rm dim}\, U_n, n_V= {\rm dim} V_n$ for the reduced spaces. This indicates the expected fact that the larger the diffusion the
smoother is the dependence of $u(\mup)$ on the parameter $\mup$. In fact, when $\e \to 0$ one approaches the
regime of pure transport where the smoothness of the parameter dependence is merely H\"older continuity requiring for a given target accuracy
a larger number of reduced basis functions, see \cite{DPW}.

%%%%%%%%%%%%%%%%%%%%%%%%%%%%%
\section{Is it Necessary to Resolve All of $\cM$?}\label{sec:duality}
%%%%%%%%%%%%%%%%%%%%%%%%%%%%%
The central focus of the preceding discussion has been to control the maximal deviation
\begin{equation}
\label{maxdist}
\sigma_n(\cM)_U = \max_{\mup\in\cY}\|u(\mup)-P_{U,U_n}u(\mup)\|_U,
\end{equation}
and to push this deviation below a given tolerance for $n$ as small as possible. However, in many applications one is not interested
in the whole solution field but only in a {\em quantity of interest} $I(\mup)$, typically of the form $I(\mup)=\ell(u(\mup))$
where $\ell\in U'$ is a bounded linear functional. Looking then for some desired optimal state $I^*= \ell(u(\mup^*))$ one is
interested in a guarantee of the form
\begin{equation}
\label{guarantee}
|\ell(u_n(\mup)) - \ell(u(\mup))|\le {\rm tol},\quad \mup \in \cY,
\end{equation}
where the states $u_n(\mup)$ belong to a possibly small reduced space $U_n$ in order to be then able to carry out the
optimization over $\mup\in \cY$ in the small space $U_n\subset U$. Asking only for the values of just a {\em linear
functional} of the solution seems to be much less demanding than asking for the whole solution and one wonders whether this can be
exploited in favor of even better online efficiency.

Trying to reduce computational complexity by exploiting the fact that, retrieving only a linear functional of an unknown state -
a scalar quantity - may require less information than recovering the whole state, is the central theme
of {\em goal oriented} adaptation in finite element methods, see \cite{BR}.
Often the desired accuracy is indeed observed to be reached by significantly coarser discretizations than needed to
approximate the whole solution within a corresponding accuracy. The underlying effect, sometimes referred to as
``squared accuracy'' is well understood and exploited in the RBM context as well, see \cite{Grepl1,Prud}.
We briefly sketch the main ideas for the current larger scope of problems and point out that, nevertheless, a guarantee of the form \eqref{guarantee} ultimately
requires controlling the maximal deviation of a reduced space in the sense of \eqref{maxdist}. Hence, an optimal sampling
of a solution manifold remains crucial.

First, a trivial estimate gives for $\ell\in U'$
\begin{equation}
\label{trivial}
|\ell(u_n(\mup)) - \ell(u(\mup))|\le \|\ell\|_{U'} \|u_n(\mup)-u(\mup)\|_U
\end{equation}
so that a control of $\sigma_n(\cM)_U$ would indeed yield a guarantee. However, the $n$ needed to drive $\|\ell\|_{U'}
\sigma_n(\cM)_U$ below ${\rm tol}$ is usually larger than necessary.

To explain the principle of improving on \eqref{trivial} we consider again a variational problem of the form \eqref{b} (suppressing any
parameter dependence for a moment)
for a pair of spaces $U,V$ where we assume now that $\kappa_{U,V'}(B)\le C_b/c_b$ is already small, possibly after
renorming an initial less favorable formulation through \eqref{hatU} or \eqref{hatV}. Let $u\in U$ again denote the exact solution of \eqref{b}.
Given a $\ell\in U'$ we wish to approximate $\ell(u)$, using an {\em approximate solution} $\bar u\in W\subset U$ defined by
\begin{equation}
\label{appru}
b(\bar u,v)= \langle f, v\rangle,\quad v\in \tilde{V}(W)\subset V,
\end{equation}
where $\tilde{V}(W)$ is a suitable test space generated by the methods discussed in \S \ref{ssec:stabvar}. In addition we will use
the solution $z\in V$ of the {\em dual problem}:
\begin{equation}
\label{dual}
b(w,z)= -\ell(w),\quad w\in U,
\end{equation}
together with an approximation $\bar z \in Z\subset V$ defined by
\begin{equation}
\label{apprdual}
b(w,\bar z)= -\ell(w),\quad w\in \tilde{W}(Z)\subset U,
\end{equation}
again with a suitable test space $\tilde{W}(Z)$. Recall that we need not determine the test spaces
$\tilde{V}(W), \tilde{W}(Z)$ explicitly but rather realize the corresponding Petrov-Galerkin projections through
the equivalent saddle point formulations with suitable $\delta$-proximal auxiliary spaces generated by
a greedy stabilization.

Then, defining the primal residual functional
\begin{equation}
\label{primres}
r_{\bar u}(v):= r(\bar u, v):= b(u-\bar u,v) = \langle f,v\rangle - b(\bar u,v),
\end{equation}
and adapting the ideas in \cite{Grepl1,Prud} for the symmetric case $V=U$ to the present slightly more general setting, we claim that
\begin{equation}
\label{hatl}
\hat\ell(\bar u):= \ell(\bar u)- r(\bar u, \bar z)
\end{equation}
is an approximation to the true value $\ell(u)$ satisfying
\begin{equation}
\label{quadratic}
|\hat\ell(\bar u)-\ell(u)| \le C \inf_{w\in W}\|u- w\|_U \inf_{v\in Z}\|z - v\|_V,
\end{equation}
where $C$ depends only on the inf-sup constant of the finite dimensional problems.
In fact, since by \eqref{dual},
\begin{equation*}
\ell(u)-\ell(\bar u)= b(\bar u-u,z) = -r(\bar u,z),
\end{equation*}
one has $\ell(u)= \ell(\bar u) - r(\bar u,z)$ and hence
\begin{eqnarray*}
|\hat\ell(\bar u)- \ell(u)| &=&| \ell(\bar u) - r(\bar u,\bar z) -\ell(\bar u)+ r(\bar u, z)|
 =| r(\bar u, z- \bar z)| = |b(u-\bar u, z-\bar z)|\\
 &\le &C_b\|u-\bar u\|_U \|z- \bar z\|_V ,
\end{eqnarray*}
which confirms the claim since $\bar u$, $\bar z$ are near-best approximations due to the asserted
inf-sup stability of finite dimensional problems.

Clearly, \eqref{quadratic} says that in order to approximate $\ell(u)$ the primal approximation in $U$ need not
resolve $u$ at all as long as the dual solution $z$ is approximated well enough. Moreover, when $\ell$ is a local
functional, e.g. a local average approximating a point evaluation, $z$ is close to the corresponding Green's function
with (near) singularity in the support of $\ell$. In the elliptic case $z$ would be very smooth away from the support
of $\ell$ and hence well approximable by a relatively small number of degrees of freedom concentrated around the
support of $\ell$. Thus, it may very well be more profitable to spend less effort on approximating $u$ than on approximating
$z$.

Returning to parameter dependent problems \eqref{parb}, the methods in \S \ref{sec:d-greedy} can now be used
as follows to construct possibly small reduced spaces for a frequent online evaluation of the quantities $I(\mup)=\ell(u(\mup))$.
We assume that we already have properly renormed families of norms $\|\cdot\|_{U_\mup}, \|\cdot\|_{V_\mup}$, $\mup\in \cY$,
with uniform inf-sup constants close to one. We also assume now that both families of norms are equivalent (by compactness of $\cY$ uniformly equivalent) to reference norms $\|\cdot\|_U, \|\cdot\|_V$, respectively. Hence we can consider two solution manifolds
\begin{equation*}
\cM_{\rm pr}:= \{u(\mup)=B_\mup^{-1}f,\, \mup\in\cY\}\subset U,\quad \cM_{\rm dual}:= \{z(\mup):= B_{\mup}^{-*}\ell,\, \mup\in\cY\}\subset V,
\end{equation*}
and use Algorithm \ref{alg:doublegreedy}, \DoubGA~ to generate (essentially in parallel) two sequences of pairs of reduced spaces
\begin{equation*}
(U_n, V_n),\, (Z_n,W_n), \quad n\in \N.
\end{equation*}
Here $V_n\subset V, W_n\subset U$ are suitable stabilizing spaces such that for $m< n$ and for the corresponding reduced solutions
$u_m(\mup)\in U_m,z_{n-m}(\mup) \in Z_{n-m}$ the quantity
\begin{equation}
\label{Inm}
I_{n,m}(\mup):= \ell(u_m(\mup))- r(u_m(\mup),z_{n-m}(\mup))
\end{equation}
satisfies
\begin{equation}
\label{quadratic2}
|I(\mup)-I_{n,m}(\mup)| \le C \sigma_m(\cM_{\rm pr})_U \sigma_{n-m}(\cM_{\rm dual})_V,
\end{equation}
with a constant $C$ independent of $n,m$. The choice of $m<n$ determines how to distribute the computational effort
for computing the two sequences of reduced bases and their stabilizing companion spaces.
By Theorem \ref{thm:main}, one can see that whichever $n$-width rate
$d_n(\cM_{\rm pr})_U$ or $d_n(\cM_{\rm dual})_V$ decays faster one can choose $m<n$ to achieve for a total of ${\rm dim}\, U_m + {\rm dim}\,Z_{n-m}
= n$ the smallest error bound. Of course, the rates are not known and one can use the tight surrogates to bound and estimate the respective errors
very accurately. For instance, when $d_n(\cM_{\rm pr})_U\le Cn^{-\alpha}$, $d_n(\cM_{\rm dual})_V\le Cn^{-\beta}$, $m= \Big\lfloor
\Big(\frac{\alpha}{\alpha +\beta}\Big)n\Big\rfloor$ yields an optimal distribution with a bound
\begin{equation}
\label{pol}
|I(\mup)-I_{n,m}(\mup)| \le C \Big(\frac{\alpha +\beta}{\beta}\Big)^\beta \Big(\frac{\alpha +\beta}{\alpha}\Big)^\alpha n^{-(\alpha +\beta)}.
\end{equation}
In particular, when $\beta > \alpha$ the dimensions on the reduced bases for the dual problem should be somewhat larger but
essentially using the same dimensions for the primal and dual reduced spaces yields the rate $n^{-(\alpha +\beta)}$ confirming
the ``squaring'' when $\alpha = \beta$. In contrast, as soon as either one of the $n$-width rates decays exponentially
it is best to grow only the reduced spaces for the faster decay while keeping a fixed space for the other side.

%%%%%%%%%%%%%%%%%%%%%%%%%%%%%
\section{Summary}
%%%%%%%%%%%%%%%%%%%%%%%%%%%%%
We have reviewed recent developments concerning reduced basis methods with the following main focus.
Using Kolmogorov $n$-width as a benchmark for the performance of reduced basis methods in terms of minimizing the dimensions
of the reduced models for a given target accuracy, we have shown that this requires essentially to construct tight well-conditioned surrogates for
the underlying variational problem. We have explained how {\em renormation} in combination with {\em inner stabilization loops}
can be used to derive such residual based surrogates
even for problem classes not covered by conventional schemes. This includes in a fully robust way
indefinite as well as ill-conditioned (singularly perturbed)
coercive problems. Greedy strategies based on such surrogates are then shown to constitute
an optimal sampling strategy, i.e., the resulting snapshots span reduced spaces whose distances from the solution manifold decay
essentially at the same rate as the Kolmogorov $n$-widths. This means, in particular, that stability constants need not be determined
by additional typically expensive computations but can be pushed by the stabilizing inner greedy loop as close to one
as one wishes.
Finally, we have explained why the focus on uniform approximation
of the entire solution manifold is equally relevant for applications where only functionals of the parameter dependent solutions
have to be approximated.

%%%%%%%%%%%%%%%%%%%%%%%%%%%%%
%%       Bibliography
%%%%%%%%%%%%%%%%%%%%%%%%%%%%%

\end{document}